\documentclass[review,onefignum,onetabnum]{siamart190516}

\usepackage{amsmath}
\usepackage{amssymb}
\usepackage{amstext}
\usepackage{amscd}
\usepackage{arydshln}
\usepackage{color}
\usepackage{graphicx}
\usepackage{longtable}
\usepackage{calc}
\usepackage{microtype}
\usepackage{alltt,xcolor}

\usepackage{amsfonts}
\usepackage{graphicx}
\usepackage{epstopdf}
\usepackage{algorithmic}
\ifpdf
\DeclareGraphicsExtensions{.eps,.pdf,.png,.jpg}
\else
\DeclareGraphicsExtensions{.eps}
\fi

\newsiamremark{remark}{Remark}
\newsiamremark{hypothesis}{Hypothesis}
\crefname{hypothesis}{Hypothesis}{Hypotheses}
\newsiamthm{claim}{Claim}
\newtheorem{example}[theorem]{Example}

\headers{Diagonal scalings for the eigenstructure of arbitrary pencils}{F. M. Dopico, M. C. Quintana, and P. Van~Dooren}

\title{Diagonal scalings for the eigenstructure of arbitrary pencils\thanks{Submitted to the editors DATE.
		\funding{Supported by ``Ministerio de Econom\'ia, Industria y Competitividad (MINECO)" of Spain and ``Fondo
			Europeo de Desarrollo Regional (FEDER)" of EU through grants MTM2015-65798-P and MTM2017-90682-REDT
			and by the ``Proyecto financiado por la Agencia Estatal de Investigaci\'on (PID2019-106362GB-I00 / AEI / 10.13039/501100011033)"  .
			The research of M. C. Quintana was funded by the ``contrato predoctoral" BES-2016-076744 of MINECO  and by an Academy of Finland grant (Suomen Akatemian p\"{a}\"{a}t\"{o}s 331240). This work was partially developed while Paul Van Dooren held a ``Chair of Excellence UC3M - Banco de Santander'' at Universidad Carlos III de Madrid in the academic year 2019-2020.}}}

\author{Froil\'{a}n M. Dopico\thanks{Departamento de Matem\'aticas,
		Universidad Carlos III de Madrid, Avda. Universidad 30, 28911 Legan\'es, Spain
		(\email{dopico@math.uc3m.es}).}
	\and Mar\'{i}a C. Quintana\thanks{Department of Mathematics and Systems Analysis, Aalto University, Otakaari 1, Espoo, Finland (\email{maria.quintanaponce@aalto.fi}).}
	\and Paul Van Dooren\thanks{Department of Mathematical Engineering, Universit\'{e} catholique de Louvain, Avenue Georges Lema\^itre 4, B-1348 Louvain-la-Neuve, Belgium
		(\email{paul.vandooren@uclouvain.be}).}
}

\usepackage{amsopn}
\DeclareMathOperator{\diag}{diag}

\def\d{\delta}

\def\1{{\mathbf{1}}}

\definecolor{myred}{cmyk}{0.000000,1.000000,1.000000,0.1}
\definecolor{myblue}{cmyk}{1.000000,0.750000,0.000000,0.1}
\newcommand{\cblu}[1]{{\color{black}{#1}}}

\newcommand{\beq}{\begin{equation}}
\newcommand{\eeq}{\end{equation}}
\newcommand{\bac}[1]{\left[\begin{array}{#1}}
	\newcommand{\eac}{\end{array}\right]}
\newcommand{\bmx}{\begin{bmatrix}}
	\newcommand{\emx}{\end{bmatrix}}

\newcommand{\la}{\lambda}

\newcommand{\hide}[1]{}

\newcounter{parenumi}
	{\end{list}}

\newlength{\tablecolwidth}
\ifpdf
\hypersetup{
  pdftitle={Diagonal scalings for the eigenstructure of arbitrary pencils},
  pdfauthor={F. M. Dopico, M. C. Quintana, and P. Van~Dooren}
}
\fi

\begin{document}

\nolinenumbers
\maketitle

\begin{abstract}
  In this paper we show how to construct diagonal scalings for arbitrary matrix pencils $\la B-A$, in which both $A$ and $B$ are complex matrices (square or nonsquare). The goal of such diagonal scalings is to ``balance'' in some sense the row and column norms of the pencil. We see that the problem of scaling a matrix pencil is equivalent to the problem of scaling the row and column sums of a particular nonnegative matrix. However, it is known that there exist square and nonsquare nonnegative matrices that can not be scaled arbitrarily. To address this issue, we consider an approximate embedded problem, in which the corresponding nonnegative matrix is square and can always be scaled. \cblu{The new scaling methods are then based on the Sinkhorn--Knopp algorithm for scaling a square nonnegative matrix with total support to be doubly stochastic or on a variant of it. In addition, using results of U. G. Rothblum and H. Schneider (1989), we give simple sufficient conditions on the zero pattern for the existence of diagonal scalings of square nonnegative matrices to have any prescribed common vector for the row and column sums}. We illustrate numerically that the new scaling techniques for pencils improve the \cblu{accuracy} of the computation of their eigenvalues.
\end{abstract}



\section{Introduction}

The problem of scaling an entrywise nonnegative $m\times n$ matrix $A$ with diagonal transformations and prespecified vectors $r$ and $c$ for the row and column sums, respectively, consists of finding a matrix of the form $S = D_{\ell}A D_r$, where $D_{\ell}\in\mathbb{R}^{m\times m}$ and $D_{r}\in\mathbb{R}^{n\times n}$ are diagonal matrices having positive diagonal elements, and such that
\begin{equation}\label{RS_problem}
S\1_n=r\quad \text{ and }\quad \1_m^{T}S=c^{T},
\end{equation}
where $\1_i:=[1,\ldots,1]^{T}\in\mathbb{R}^{i}$ for $i=n,m$ \cite{Brualdi66,RoSc}. When $r=\1_m$ and $c=\1_n$ the scaled matrix $S$ is \cblu{neccessarily square and is} said to be doubly stochastic, i.e., its row and column sums are all equal to $1$.

The related problem of scaling the rows and columns of a complex square matrix $A$ (not necessarily nonnegative) using real and positive diagonal similarity transformations in order to compute more accurate eigenvalues, is a well established technique to improve the sensitivity of the eigenvalue problem of the matrix $A$ \cite{Par}. \cblu{This is known as {\em balancing} the matrix $A$. In exact arithmetic, it amounts to minimizing the Frobenius norm of the scaled matrix $D^{-1}AD$, where $D$ ranges over all non-singular real diagonal matrices, which is equivalent to minimizing the departure from normality of $D^{-1}AD$ \cite{LemVD}. Since the eigenvalues of normal matrices have condition numbers equal to $1$, such scaling very often improves the sensitivity of eigenvalues}. The method for computing the optimal scaling is a very simple cyclic procedure where at each step only a single diagonal element of $D$ is updated. This method is implemented in MATLAB \cite{Matlab} as a default option of the eigenvalue computation problem, which indicates that its effectiveness is well accepted.  For improving  the accuracy of the eigenvalues computed in floating point arithmetic, it is essential that the diagonal elements of $D$ are integer powers of $2$, because in this way the scaling does not produce any rounding errors and the eigenvalues are preserved {\em exactly} under such a scaling transformation. Otherwise, the rounding errors inherent to constructing $D^{-1}AD$ would spoil any potential improvement in the accuracy of the computed eigenvalues. As explained in \cite{Par}, the restriction to diagonal matrices $D$ whose entries are integer powers of $2$ allows for a relaxed stopping criterion of the cyclic procedure for computing $D$ and implies that the related minimization problem is only approximately solved.

The idea of performing positive diagonal scalings in order to improve the \cblu{accuracy of computed} eigenvalues was also extended to the generalized eigenvalue problem of a regular pencil  $\lambda B-A$. In this case, the nonsingular diagonal matrices multiplying the pencil on the left and on the right are different. In \cite{Ward}, Ward describes a scaling technique which aims at making the pencil entries have magnitudes as close to unity as possible. \cblu{In \cite{LemVD}, Lemonnier and Van Dooren propose a diagonal scaling that in exact arithmetic minimizes the Frobenius norm of the pencils over all positive diagonal scalings with fixed determinant. This improves very often the conditioning of the eigenvalues, since the solution of such minimization problem over general nonsingular transformations
is a so-called standardized normal pencil, which is a pencil whose eigenvalues all have a condition number in the chordal metric that is smaller than or equal to $\sqrt{2}$. The method of Ward is the one that LAPACK \cite{lapack} proposes as built-in option for scaling a regular pencil, but it was pointed out in \cite{LemVD} that the method of Lemonnier-Van Dooren outperforms that of Ward in terms of the accuracy of the computed eigenvalues, especially when the pencil has entries of strongly varying magnitudes. The experiments in Section \ref{sec:numerics} will further confirm the superiority of the method in \cite{LemVD} for a wide variety of pencils of different sizes and types. As in the case of balancing matrices, it is essential that the entries of the diagonal scaling matrices are integer powers of $2$ in order to improve the accuracy of the computed eigenvalues in floating point arithmetic. Currently, MATLAB does not offer any built-in option for scaling pencils. We will see in Section \ref{sec:regular} that the method in \cite{LemVD} is equivalent to scaling a particular nonnegative matrix to a multiple of a doubly stochastic matrix, which motivates us to revise briefly the literature on this and other related problems.

There is a vast literature on diagonal scaling of nonnegative matrices for getting a matrix with prescribed row and column sums. The origin of these problems goes back at least until the beginning of the XX century \cite{kruithof37,yule2012} and originates in the area of optimal transport \cite{CuturiPeyre}, though it has applications in many other areas \cite{Idel}. See \cite[Section 3.1]{Idel} and  \cite[Remark 4.5]{CuturiPeyre} for historical remarks on these problems. Relevant classical references from the point of view of matrix analysis include
\cite{Brualdi66,Krupp,RoSc,SiK}, among many others.  Despite this vast literature, several issues are still open for improvement, such as a good understanding of the convergence of related algorithms for sparse matrices and simple conditions on the zero pattern of the matrix for existence and unicity of a solution for special cases, specially in the case of rectangular matrices \cite{Idel}. The most relevant papers on diagonal scalings that are closely related to the problems discussed in this paper are, in chronological order, those of Sinkhorn-Knopp \cite{SiK}, Krupp \cite{Krupp}, Rothblum-Schneider \cite{RoSc} and Knight \cite{Knight}, which is why we quote theorems from those papers.

In this paper we show that there exists a link between the problem of scaling a regular square pencil and that of scaling a square nonnegative matrix to become doubly stochastic. This implies that the scaling is essentially unique and bounded if and only if the corresponding nonnegative matrix satisfies certain conditions, namely {\em total support} and {\em full indecomposability}.} Moreover, in that situation, the scaling can be found through the well-known Sinkhorn-Knopp algorithm \cite{Knight,SiK}. We then show how to extend this to singular or nonsquare pencils, which, to the best of our knowledge, has not been considered yet in the literature. For that, we introduce a regularization term into the original problem which ensures existence of a solution of an approximate problem with bounded diagonal scalings $D_\ell$ and $D_r$. In addition, the regularization term can be considered in both square or nonsquare cases.

These ideas are connected to the results of Rothblum and Schneider \cite{RoSc} about \cblu{scaling arbitrary nonnegative matrices (square or rectangular) with} prespecified row and column sums,\cblu{ which can be obtained using a Sinkhorn-Knopp-like algorithm, but many other optimization methods have been proposed in the literature \cite{Idel,CuturiPeyre}.}
We then build on these ideas to further improve the scaling technique of Lemonnier and Van Dooren by introducing the regularization term as an additional cost. This cost can be viewed as a regularization to ensure \cblu{always the} existence and boundedness of our scaling, but it also ensures essential unicity of the computed scaling.

The paper is organized as follows. In Section \ref{sec:scaling}, we give some basic notions about scaling pencils\cblu{, scaling nonnegative matrices and the Sinkhorn-Knopp-like algorithm}. In Sections \ref{sec:regular} and \ref{sec:nonsquare}, we study the diagonal scaling problem for square and nonsquare pencils, respectively. In Section \ref{sec:regular}, we will also recall the necessary and sufficient conditions for a square \cblu{nonnegative} matrix to become doubly stochastic under diagonal scalings, and we give \cblu{simple sufficient conditions based on the zero pattern of the matrix} for the existence of diagonal scalings having any prespecified common vector for the row and column sums. These results will be useful in Section \ref{sec:regularized}.  In that section, we develop a new scaling technique for generalized eigenvalue problems and show that it can be applied to any pencil, regular or singular, square or rectangular. For that, we introduce a regularization term into the original problem which guarantees existence, unicity and boundedness of the scaling. In addition, in Subsection \ref{sec:regularized_rowandcolum}, we consider a modified version of the new scaling technique that is \cblu{often} better for scaling nonsquare pencils. In Section \ref{sec:numerics} we then illustrate the improved accuracy of the computed eigenvalues using several numerical examples. In the last Section \ref{sec:conclusion} we give some concluding remarks.

\section{Preliminaries: \cblu{Scaling arbitrary pencils and nonnegative matrices}}  \label{sec:scaling}

The standard techniques for computing eigenvalues of complex pencils of matrices guarantee that the backward errors corresponding
to the computed spectrum are essentially bounded by the norm of the coefficients of the pencil, times the machine precision of the computer used.
But one can improve this bound by reducing the norms of the coefficients without affecting the spectrum. This is where balancing using diagonal scaling comes in. \cblu{We emphasize again that the diagonal entries of such scalings must be integer powers of $2$, since otherwise the rounding errors of floating point arithmetic would destroy any potential improvement in accuracy that such scalings might achieve.}

Two types of scalings can be applied to a pencil $\la B-A$.

The first one is a change of variable $\hat \la := d_{\la}\la$ to make sure that the scaled matrices $A$ and $B/d_{\la}$ have approximately the same norm. This can be done without introducing rounding errors, by taking $d_\la$ equal to a power of 2. The staircase and the $QZ$ algorithm work independently on both matrices and this scaling can be restored afterwards, again without introducing any additional errors. One could therefore argue that this scaling is irrelevant for these algorithms, but we will see that it affects the second scaling procedure we will discuss. Therefore we will assume in the sequel that both matrices $A$ and $B$ are of comparable norms, and that no such variable scaling needs to be applied.

The second type of scaling is based on multiplication on the left and on the right by positive diagonal matrices $D_\ell$ and $D_r$, respectively, that are chosen to ``balance'' in some sense the row and column norms of the complex matrices $\widetilde{A}:=D_\ell A D_r$ and $\widetilde{B}:=D_\ell B D_r$. We will see that balancing the row and column norms of the matrices $\widetilde{A}$ and $\widetilde{B}$ is equivalent to performing two-sided diagonal scalings to a particular real entrywise nonnegative matrix $M$. Therefore, we recall in the sequel some results on this problem.

The first result we revise appears in \cite[Theorem 2, (a)-(b)]{RoSc} and is the next one.
\begin{theorem} \label{thm.RS2ab}
Given a real nonnegative matrix $M\in \mathbb{R}^{m\times n}$ and vectors $r\in \mathbb{R}^{m\times 1}$ and $c\in \mathbb{R}^{n\times 1}$ with strictly positive entries satisfying $\1_m^Tr=c^T\1_n$, there exist positive diagonal matrices $D_{M,\ell}$ and $D_{M,r}$ such that
	\begin{equation}\label{eq:RS_prob}
D_{M,\ell}  M D_{M,r} \1_n  = r\quad  \text{ and }\quad  \1_m^T D_{M,\ell}  M D_{M,r}  = c^T
	\end{equation}
if and only if there exists a matrix $S$ with the same zero pattern as $M$ such that $S\1_n=r$ and $\1_m^TS=c^T$.
\end{theorem}
This is an elegant nontrivial existence result that in a less general form appeared before in \cite{menon1968}. To tackle the problem of finding the scaled matrix, one can perform a {\em Sinkhorn-Knopp-like algorithm} by alternatively normalizing the row and column sums of $M$ as follows:
\medskip

\noindent
{\bf Algorithm 1} {\em (Sinkhorn-Knopp-like algorithm for nonnegative $M \in \mathbb{R}^{m\times n}$)}

\noindent
Initialize: $D_{M,\ell} = I_m$ and $D_{M,r} = I_n$
\begin{itemize}
	\item[\rm(1)] Multiply each row $i$ of $M$ and of $D_{M,\ell}$ by $\dfrac{r_i}{\sum_{j} m_{ij}}$ to obtain an updated matrix $M$ with row sums $r$ and an updated matrix $D_{M,\ell}$.
	\item[\rm(2)] Multiply each column $j$ of the updated $M$ and of $D_{M,r}$ by $\dfrac{c_j}{\sum_{i} m_{ij}}$ to obtain an updated matrix $M$ with column sums $c$ and an updated matrix $D_{M,r}$.
	\item[\rm(3)] If the row sums of the matrix $M$ obtained in step $\rm(2)$ are far from $r$, repeat steps $\rm(1)$ and $\rm(2)$ with such $M$ until an adequate stopping criterion is satisfied.
\end{itemize}

\medskip
\noindent
We give a MATLAB code of this algorithm in Appendix $A$. This algorithm appeared as early as in \cite{yule2012} and \cite{kruithof37} and, according to \cite[Section 3.1]{Idel}, it has been rediscovered several times in the literature and has received different names as, for instance, the Kruithof's projection method (see \cite{Krupp}) or the RAS method, among many others. In this paper, we have decided to refer to this method as the Sinkhorn-Knopp-like algorithm, because if $r=c= \1_n$ and $M$ is square, then it collapses to the famous Sinkhorn-Knopp algorithm for scaling a nonnegative matrix to a doubly stochastic matrix \cite{SiK}. If the Sinkhorn-Knopp-like algorithm converges, i.e., $M$ converges and the diagonal matrices of the iteration converge to positive bounded diagonal matrices, the limit will be the scaled matrix $D_{M,\ell}  M D_{M,r}$ in Theorem \ref{thm.RS2ab}.

\cblu{Another important result in this context is that there exists at most one solution for the two-sided diagonal scaling problem in \eqref{eq:RS_prob} for any prescribed vectors $r$ and $c$. This is stated in the following Theorem \ref{RS}, which is a partial result of what is proven in \cite[Theorem 4]{RoSc}.}

\begin{theorem} \label{RS} Let $M\in \mathbb{R}^{m\times n}$ be a nonnegative matrix and let $r\in \mathbb{R}^{m\times 1}$ and $c\in \mathbb{R}^{n\times 1}$ be strictly positive vectors satisfying $\1_m^Tr=c^T\1_n$. Then there exists at most one two-sided scaled matrix $S= D_{M,\ell} M D_{M,r}$ with row sums $S\1_n=r$ and column sums $\1_m^TS=c^T$, where $D_{M,\ell}$ and $D_{M,r}$ are diagonal matrices with positive main diagonals.
\end{theorem}
\cblu{A less general version of Theorem \ref{RS} appeared in \cite{menon1968} and the general case is implicit in \cite{menon1969}. We emphasize that, although $S$ is unique when it exists, the matrices $D_{M,\ell}$ and $D_{M,r}$ are not necessarily unique. We refer the reader to \cite[Theorem 4]{RoSc} for a description of all matrices $D_{M,\ell}$ and $D_{M,r}$ that satisfy $S= D_{M,\ell} M D_{M,r}$.}

A surprising and useful result is that the Sinkhorn-Knopp-like algorithm converges if and only if the scaling problem \eqref{eq:RS_prob} has solution. This was proved for general matrices and arbitrary prescribed row and column sum vectors in \cite{Krupp} and for the case of square nonnegative matrices and $r = c = 1_n$ in \cite{SiK}, i.e., for the doubly stochastic case (see also \cite[Theorem 4.1]{Idel}). Next, we state this important result.

\begin{theorem}\label{conv} Under the assumptions in Theorem \ref{RS}, there exist diagonal matrices $D_{M,\ell}$ and $D_{M,r}$ with positive main diagonals such that \eqref{eq:RS_prob} is satisfied if and only if the Sinkhorn-Knopp-like algorithm converges.
\end{theorem}

Therefore, if a nonnegative matrix $M$ can be scaled for prescribed row and column sums, the scaled matrix is unique  and is the limit of the Sinkhorn-Knopp-like algorithm, which gives a practical numerical procedure to check for scalability. Unfortunately, the Sinkhorn-Knopp-like algorithm can be very slow, in particular for sparse matrices, and other faster algorithms have been developed in the literature (see \cite[Section 7]{Idel}, \cite[Section 4.3]{CuturiPeyre} and the references therein). However, we emphasize that for the main purpose of this paper, i.e., improving the accuracy of computed eigenvalues of pencils, we have always found that the Sinkhorn-Knopp-like algorithm is fast enough and that the cost of its application is much smaller than the cost of computing the eigenvalues. The reason is that, in this case, the diagonal entries of the scalings $D_{M,\ell}$ and $D_{M,r}$ to be applied to the pencil must be integer powers of $2$ which allows to use a very relaxed stopping criterion in the Sinkhorn-Knopp-like algorithm. We will discuss this issue in depth in Section \ref{sec:numerics}.

One can find necessary and sufficient non-algorithmic conditions for the scaled matrix to exist in \cite[Theorem 2]{RoSc}, \cite[Theorem 2.1]{Brualdi66} and \cite[Theorem 4.1]{Idel}. However, these conditions depend on nontrivial properties that must be satisfied by the vectors $r$ and $c$, as those we state in Lemma \ref{lem_RS}. In general, necessary and sufficient conditions depending only on the zero pattern of $M$ are not known. A remarkable exception to this comment is the doubly stochastic scaling problem $r = c = \1_n$ for square matrices, where such a condition is provided by the total support of the matrix (see Section \ref{sec:regular}). In the next section, we will present new simple sufficient conditions depending only on the zero pattern for diagonal scalings to exist with prescribed common vector for the row and column sums in the case of balancing square pencils and matrices.

There are infinitely many examples of nonnegative matrices that cannot be scaled for prescribed $r$ and $c$. \cblu{The following example illustrates this fact.}

\begin{example} \label{ex:nonsquare} For instance, one can easily check that the matrix
	$$M:=\begin{bmatrix}
	1 & 1 & 1 \\
	0 & 0 & 1
	\end{bmatrix}$$ can not be scaled with prescribed vectors $r:=[3, 3]^{T}$, for the row sums, and  $c:=[2, 2, 2]^{T}$, for the column sums.
\end{example}

\section{Scaling square pencils and related problems} \label{sec:regular}

Let us first look at the case of square pencils.
In \cite[page 259]{LemVD}, positive diagonal matrices $D_\ell$ and $D_r$ are chosen to equilibrate the row and column norms of a $n\times n$ regular pencil $\la B-A$, by imposing
\begin{equation} \label{balanced}
\|\text{col}_j (\widetilde{A})\|_2^2 + \|\text{col}_j (\widetilde{B})\|_2^2 = \|\text{row}_i (\widetilde{A})\|_2^2 + \|\text{row}_i (\widetilde{B})\|_2^2 = \gamma^2 \text{, for }i,j=1,\ldots,n,
\end{equation}
for some constant $\gamma$ resulting from the balancing, where $\widetilde{A}:=D_\ell A D_r$ and $\widetilde{B}:=D_\ell B D_r$, and $\|\cdot\|_2$ denotes the standard Euclidean norm of a vector \cite{GoVL}.
A pencil satisfying these conditions was called {\em balanced} and an algorithm was presented in \cite{LemVD} to compute a scaling to balance a regular pencil
$\lambda B-A$. It was shown that this amounts to solving the following norm minimization problem
\begin{equation} \label{inf1}
\inf_{\det{D_\ell}.\det{D_r}=1} \|D_\ell(\la B-A)D_r\|_F^2,
\end{equation}
using the so-called Frobenius norm of a pencil:
$$  \|\la B-A\|_F^2 :=  \|B\|_F^2+ \|A\|_F^2,
$$
where $ \|A\|_F$ and $\|B\|_F$ are the matrix Frobenius norms of $A$ and $B$ \cite{GoVL}. Moreover, the following result was proven in \cite{LemVD}.
\begin{theorem} The minimization problem
	\begin{equation} \label{infT}
	\inf_{\det{T_\ell}.\det{T_r}=1} \|T_\ell(\la B-A)T_r\|_F^2,
	\end{equation}
	where $T_\ell$ and $T_r$ are arbitrary nonsingular matrices, has a so-called {\em standardized normal pencil} $\la \hat B- \hat A$ as solution, satisfying
	$$ U_\ell(\la \hat B- \hat A)U_r =  \la \Lambda_B-\Lambda_A, \quad U_\ell^* U_\ell= U_r^* U_r=I_n, \quad |\Lambda_B|^2+ |\Lambda_A|^2 = \gamma^2 I_n,
	$$
	where $\Lambda_B$ and $\Lambda_A$ are diagonal. If the eigenvalues of the regular pencil $\la B-A$ are distinct, then $T_\ell$ and $T_r$ have a bounded solution and the infimum is a minimum; otherwise they may be unbounded.
\end{theorem}
As shown in \cite{LemVD}, the standardized normal pencils happen to have eigenvalues with condition number bounded by $\sqrt{2}$.
This explains why performing the same minimization over the diagonal scalings is likely to improve the sensitivity of the eigenvalue computation. Moreover, if the transformation matrices are bounded then the eigenstructure of the regular pencil is preserved.

\medskip

But the positive diagonal scalings that achieve the balancing in \cite{LemVD} are not unique, and they may not exist or may be unbounded.
In order to analyze this further we relate this problem to that of scaling a real nonnegative square matrix by two-sided scalings to a doubly stochastic matrix, or in other words, to make the row sums and column sums equal to 1. \cblu{As mentioned before, an algorithm to solve this problem has been developed and
analyzed by Sinkhorn and Knopp \cite{SiK} and reduces to Algorithm 1 with $r = c = \1_n$. Further analysis can be found in \cite{Knight}}. The link between both problems is the following.
Let us define the nonnegative matrices
\begin{equation} \label{nonneg}
M :=  |A|^{\circ 2}+|B|^{\circ 2}, \quad \mathrm{and} \quad \widetilde{M} :=  |\widetilde{A}|^{\circ 2}+|\widetilde{B}|^{\circ 2}
\end{equation}
where $|X|$ indicates the element-wise absolute value of the matrix $X$, where $X^{\circ 2}$ indicates the elementwise square of the matrix $X$,
and where $D_\ell$ and $D_r$ satisfy the balancing equations \eqref{balanced}.
Then the scaled matrix  $ \widetilde{M}= D_\ell^2 M  D_r^2 $ satisfies
$$   \widetilde{M} \mathbf{1}_n= D_\ell^2 (|A|^{\circ 2}+|B|^{\circ 2})  D_r^2 \mathbf{1}_n= \gamma^2\mathbf{1}_n , \quad
\mathbf{1}_n^T\widetilde{M}= \mathbf{1}_n^T D_\ell^2 (|A|^{\circ 2}+|B|^{\circ 2})  D_r^2 =  \gamma^2\mathbf{1}_n^T
$$
which implies that $\widetilde{M}/\gamma^2$ is doubly stochastic and that the two-sided scaling for the nonnegative matrix $M$ satisfies
$$ \widetilde M /\gamma^2= D_{M,\ell} M D_{M,r}, \quad \mathrm{where} \quad
D_{M,\ell}:=D_\ell^{2}/\gamma, \; D_{M,r}:=D_r^2/\gamma.$$
The only difference is that for balancing, we impose a scalar constraint $\det{D_\ell} \cdot \det{D_r}=1$, which is why the resulting row and column norms are equal to $\gamma^2$ rather than 1. In fact, the algorithm proposed in \cite{LemVD} was to alternately normalizing the rows and columns of $M$ to 1 (rather than $\gamma$), and that is precisely the algorithm of Sinkhorn-Knopp. This connection was not established in \cite{LemVD}.

It follows from this that the unicity or boundedness of the scalings are equivalent for the two problems.

We recall in Theorem \ref{SK} the results given for two-sided scaling in \cite{SiK} for square nonnegative matrices $M\in \mathbb{R}^{n\times n}$ in order \cblu{for} the corresponding matrix to become doubly stochastic. We notice that the doubly stochastic scaling problem of Theorem \ref{SK} is a special case of the scaling problem in Theorem \ref{RS}, just by considering square matrices and $r=c=\1_n$. Before stating Theorem \ref{SK}, we introduce the notions of total support and full indecomposability, that will be used.

\begin{definition}
	The sequence $m_{1,\sigma(1)},m_{2,\sigma(2)},\cdots,m_{n,\sigma(n)}$, where  $\sigma$ is a permutation of $\{1,2,\cdots,n\}$,
	is called a diagonal of a $n\times n$ square matrix $M$. A nonnegative matrix $M\in\mathbb{R}^{n\times  n}$ is said to have {\em total support} if every positive element of $M$ lies on a positive diagonal.
\end{definition}

\begin{definition}
	A nonnegative matrix $M\in\mathbb{R}^{n \times  n}$ is said to be {\em fully indecomposable} if there do not exist permutation matrices $P_\ell$ and $P_r$ such that
	$P_\ell M P_r$ can be partitioned as
	$$  P_\ell M P_r=\left[\begin{array}{cc} M_{11} & M_{12} \\ 0 & M_{22} \end{array} \right],
	$$
	where $M_{11}$ and $M_{22}$ are square matrices.
\end{definition}
\begin{remark} \rm It was proved in \cite{Brualdi2} that a fully indecomposable matrix has total support.
\end{remark}

\begin{theorem}(Sinkhorn-Knopp) \label{SK}
	If $M\in \mathbb{R}^{n\times n}$ is a nonnegative matrix then a necessary and sufficient condition that there exists a doubly stochastic matrix $S$ of the form $S= D_{M,\ell} M D_{M,r}$, where $D_{M,\ell}$ and $D_{M,r}$ are diagonal matrices with positive main diagonals, is that $M$ has total support. If $S$ exists, then it is unique. $D_{M,\ell}$ and $D_{M,r}$ are also unique up to a nonnegative scalar multiple if and only if $M$ is fully indecomposable.
\end{theorem}

The doubly stochastic matrix $S$ can be obtained as a limit of a sequence of matrices generated by alternately normalizing the row and column sums of $M$\cblu{, i.e., by applying Algorithm 1 with $r = c = \1_n$, which is the Sinkhorn-Knopp algorithm. As a consequence of Theorems \ref{conv} and \ref{SK}, a necessary and sufficient condition that the Sinkhorn-Knopp algorithm applied to $M$} will converge to a doubly stochastic limit of the form $D_{M,\ell} M D_{M,r}$ is that $M$ has total support \cite{Knight,SiK}. 

We recall in the following Theorem \ref{th:symmetric} the particular case of having a symmetric and fully indecomposable matrix $M$. This case will be important in the new regularized scaling method developed in Section \ref{sec:regularized}.

\begin{theorem}\cite[Lemma 4.1]{Knight}\label{th:symmetric} If $M\in \mathbb{R}^{n\times n}$ is a symmetric nonnegative and fully indecomposable matrix then there exists a unique diagonal matrix $D$ with positive main diagonal such that $DMD $ is doubly stochastic.
\end{theorem}

\begin{remark}\label{bounded} When $M$ is fully indecomposable, the solution set for the diagonal scalings is $\mathcal{S}:=\{(D_{M,\ell}/c, cD_{M,r}) : c > 0\}$, for a given solution $(D_{M,\ell}, D_{M,r})$. To guarantee unicity for a solution in $\mathcal{S}$, one can consider a unique ``normalized'' scaling pair $(D_{M,\ell}, D_{M,r})$. For instance, by imposing that the solution satisfies $\det D_{M,\ell} = \det D_{M,r}$ or $\displaystyle\max_{i=1,\ldots,n} \{d_{i}^{\ell}\} = \displaystyle\max_{i=1,\ldots,n}\{d_{i}^{r}\} $, where $d_{i}^{\ell}$ and $d_{i}^{r}$ are the diagonal entries of $D_{M,\ell}$ and $D_{M,r}$, respectively. Then the pair $(D_{M,\ell}, D_{M,r})$ is unique in $\mathcal{S}$. Moreover, when $M$ is symmetric, then these normalizations imply that $D_{M,\ell}=D_{M,r}$. In summary, one can always perform a normalization in order to obtain unicity for the diagonal scalings. \end{remark}

In the following examples, we illustrate what is happening when the conditions mentioned in Theorem \ref{SK} do not hold.

\begin{example}\label{examples} Let us consider the regular pencil
	$$  \la B_1 - A_1 :=\left[\begin{array}{ccc} 1 & \la & 0 \\
	\la & 0 & 0 \\
	0 & 0 & 1
	\end{array}\right],\quad  \text{and let}\quad  M_1 :=\left[\begin{array}{ccc} 1 & 1 & 0 \\
	1 & 0 & 0 \\
	0 & 0 & 1
	\end{array}\right]
	$$
	be the corresponding matrix $M:=M_1$ in \eqref{nonneg}. $M_1$ has no total support since the (1,1) entry is not on a positive diagonal. The Sinkhorn-Knopp algorithm does not converge for this example. In fact, any candidate pair of scalings  $D_{M,\ell}=\diag (\ell_1,\ell_2,\ell_3)$, and
	$D_{M,r}=\diag (r_1,r_2,r_3),$ has to satisfy $\ell_1 r_2=\ell_2 r_1=\ell_3 r_3=1$ and $\ell_1 r_1=0$ which does not have a bounded solution.
	
	Now, let us consider the regular pencil
	$$  \la B_2 - A_2 :=\left[\begin{array}{ccc} 1 & \la & 0 \\
	\la & 1 & 0 \\
	0 & 0 & 1
	\end{array}\right],\quad  \text{ and let}\quad  M_2 :=\left[\begin{array}{ccc} 1 & 1 & 0 \\
	1 & 1 & 0 \\
	0 & 0 & 1
	\end{array}\right]
	$$
	be the corresponding matrix $M:=M_2$ in \eqref{nonneg}. In this case, $M_2$ has total support and the Sinkhorn-Knopp algorithm converges. Indeed, the following positive diagonal scaling makes $M$ doubly stochastic: $$ \left[\begin{array}{ccc}
	\sqrt{\frac12} & 0 & 0 \\
	0 & \sqrt{\frac12} & 0 \\
	0 & 0 & 1
	\end{array}\right] \left[\begin{array}{ccc} 1 & 1 & 0 \\
	1 & 1 & 0 \\
	0 & 0 & 1
	\end{array}\right] \left[\begin{array}{ccc}
	\sqrt{\frac12} & 0 & 0 \\
	0 & \sqrt{	\frac12}& 0 \\
	0 & 0 & 1
	\end{array}\right]  = \left[\begin{array}{ccc} \frac12 & \frac12 & 0 \phantom{\Big|}  \\
	\frac12 & \frac12 & 0 \phantom{\Big|} \\
	0 & 0 & 1
	\end{array}\right]. $$
	However, $M_2$ is not fully indecomposable, which implies that $D_{M,\ell}$ and $D_{M,r}$ are not unique up to a scalar multiple. In this case, the Sinkhorn-Knopp algorithm may converge to different diagonal scaling matrices for different starting diagonal initial conditions. Moreover, it may converge to unbounded $D_{M,\ell}$ and $D_{M,r}$.
	For instance, for the following scaling
	$$ \left[\begin{array}{ccc} t \sqrt{\frac12} & 0 & 0 \\
	0 & t \sqrt{\frac12} & 0 \\
	0 & 0 & 1/s
	\end{array}\right] \left[\begin{array}{ccc} 1 & 1 & 0 \\
	1 & 1 & 0 \\
	0 & 0 & 1
	\end{array}\right] \left[\begin{array}{ccc}
	\dfrac{1}{t}\sqrt{\frac12} & 0 & 0 \\
	0 & 	\dfrac{1}{t}\sqrt{\frac12}& 0 \\
	0 & 0 & s
	\end{array}\right]  = \left[\begin{array}{ccc} \frac12 & \frac12 & 0 \phantom{\Big|}  \\
	\frac12 & \frac12 & 0 \phantom{\Big|} \\
	0 & 0 & 1
	\end{array}\right] $$
	the right diagonal matrix is unbounded as $t\to 0$ and the left one as $s\to 0$. Finally, let us consider the regular pencil
	$$  \la B_3 - A_3 :=\left[\begin{array}{ccc} 1 & \la & 0 \\
	\la & 0 & \la \\
	0 & \la & 1
	\end{array}\right],\quad  \text{ and let}\quad  M_3 :=\left[\begin{array}{ccc} 1 & 1 & 0 \\
	1 & 0 & 1 \\
	0 & 1 & 1
	\end{array}\right]
	$$
	be the corresponding matrix $M:=M_3$ in \eqref{nonneg}. In this case, $M_3$ has total support and is, in addition, fully indecomposable. Then the scaling procedure converges to bounded diagonal scaling matrices, that are essentially unique (up to a scalar multiple):
	$$ \left[\begin{array}{ccc}
	\sqrt{\frac12} & 0 & 0 \\
	0 & \sqrt{\frac12} & 0 \\
	0 & 0 & \sqrt{\frac12}
	\end{array}\right] \left[\begin{array}{ccc} 1 & 1 & 0 \\
	1 & 0 & 1 \\
	0 & 1 & 1
	\end{array}\right] \left[\begin{array}{ccc}
	\sqrt{\frac12} & 0 & 0 \\
	0 & \sqrt{\frac12}& 0 \\
	0 & 0 & \sqrt{\frac12}
	\end{array}\right]  = \left[\begin{array}{ccc} \frac12 & \frac12 & 0 \phantom{\Big|}  \\
	\frac12 & 0 & \frac12 \phantom{\Big|} \\
	0 & \frac12 & \frac12  \phantom{\Big|}
	\end{array}\right]. $$	
\end{example}

For the general scaling problem in Theorem \ref{RS}, with arbitrary prespecified vectors for the row and column sums, sufficient conditions on $M$ for the scaling to exist \cblu{as simple as those in Theorem \ref{SK}, which are based only on the zero pattern of $M$,} are not known in the literature, \cblu{to the best of our knowledge}, not even in the case of a square matrix $M$. \cblu{This motivated us to develop the results in the next subsection.}

\subsection{\cblu{Diagonal scalings of square nonnegative matrices with prescribed common vector for the row and column sums}}
We now derive \cblu{simple} sufficient conditions \cblu{on the zero pattern} for the existence of a diagonal scaling of a square matrix $M$ by considering not only the vector $\1_n$ but any prescribed common vector \cblu{$v$} for the row and column sums. For that, we use the following Lemma \ref{lem_RS}, which is a partial result of \cite[Theorem 2]{RoSc}. In what follows, the support of a matrix $A\in \mathbb{R}^{m\times n}$, denoted by $\text{supp}(A)$, is defined as the set $\{(i, j)\,|\,a_{ij}\neq 0, i=1,\cdots,m, \text{ and } j=1,\cdots,n\}$.

\begin{lemma}\label{lem_RS} Let $M\in \mathbb{R}^{m\times n}$ be a nonnegative matrix and let $r\in \mathbb{R}^{m\times 1}$ and $c\in \mathbb{R}^{n\times 1}$ be strictly positive vectors such that $\1_m^Tr=c^T\1_n$. Then there exists a scaled matrix $S= D_{M,\ell} M D_{M,r}$ with row sums $S\1_n=r$ and column sums $\1_m^TS=c^T$, where $D_{M,\ell}$ and $D_{M,r}$ are diagonal matrices with positive main diagonals, if and only if there exist no pair of vectors $(u, v)\in \mathbb{R}^{m} \times \mathbb{R}^{n}$ for which
	\begin{itemize}
		\item[\rm(a)] $u_i + v_j \leq 0$ for each pair $(i,j)\in\text{supp}(M),$
		\item[\rm(b)] $ r^T u = c^T v = 0,$ and
		\item[\rm(c)] $u_{i_0} + v_{j_0} < 0$ for some pair $(i_0,j_0)\in\text{supp}(M)$.
	\end{itemize}		
\end{lemma}

\begin{theorem}\label{th_squareRS}   Let $M\in \mathbb{R}^{n\times n}$ be a nonnegative matrix with $(i,i)\in\text{supp}(M)$ for all $i=1,\cdots,n$ and such that $\text{supp}(M)=\text{supp}(M^T) $. Let $v\in \mathbb{R}^{n\times 1}$ be a strictly positive vector. Then there exists a scaled matrix $S= D_{M,\ell} M D_{M,r}$ with row sums $S\1_n=v$ and column sums $\1_n^TS=v^T$, where $D_{M,\ell}$ and $D_{M,r}$ are diagonal matrices with positive main diagonals. Moreover, $S$ is unique  and is the limit of the Sinkhorn-Knopp-like algorithm. If, in addition, $M$ is fully indecomposable then $D_{M,\ell}$ and $D_{M,r}$ are also unique up to a nonnegative scalar multiple \cblu{and, if $M=M^T$, then there exists a unique diagonal matrix $D$ with positive diagonal entries such that $S = DMD$}.
\end{theorem}
\begin{proof} Consider a $n\times n$ nonnegative matrix $M$ such that $\text{supp}(M)=\text{supp}(M^T) $ and $(i,i)\in\text{supp}(M)$ for all $i=1,\cdots,n$.
	By contradiction, let us assume that there exists no scaled matrix $S$ with row sums $S\1_n=v$ and column sums $\1_n^TS=v^T$. Then, by Lemma \ref{lem_RS}, there exists a pair of vectors $(x, y)\in \mathbb{R}^{n} \times \mathbb{R}^{n}$ for which
	\begin{itemize}
		\item[\rm(a)] 	$x_i + y_j \leq 0 $ for each pair $ (i,j)\in\text{supp}(M),$
		\item[\rm(b)] $	v^T x = v^T y = 0,$ and
		\item[\rm(c)] $ x_{i_0} + y_{j_0} < 0 $ for some pair $ (i_0,j_0)\in\text{supp}(M).$
	\end{itemize}
	Condition $\rm(b)$ implies that
	\begin{equation}\label{eq:sum}
	v_1(x_1+y_1)+\cdots+v_n(x_n+y_n)=0.
	\end{equation}
	In addition, since $(i,i)\in\text{supp}(M)$ for all $i=1,\dots,n$, condition $\rm(a)$ implies that $x_i+y_i\leq 0$ for all $i=1,\dots ,n$. It then follows from \eqref{eq:sum} that $x_i+y_i = 0$ for all $i=1,\dots ,n$ since $v_i>0$. Moreover, by condition $\rm(c)$, there exists a pair $ (i_0,j_0)\in\text{supp}(M)$ such that $ x_{i_0} + y_{j_0} < 0 $. Taking into account that $x_i+y_i=0$ for all $i=1,\dots , n$ we have that
	\begin{equation}\label{eq:sum2}
	(x_{i_0}+y_{i_0})+(x_{j_0}+y_{j_0})=0.
	\end{equation}
	By equation \eqref{eq:sum2} and the fact that $ x_{i_0} + y_{j_0} < 0 $, we obtain that $ x_{j_0} + y_{i_0} > 0 $. Therefore, by $\rm(a)$, $(j_0,i_0)\not\in\text{supp}(M)$, which is a contradiction since $ (i_0,j_0)\in\text{supp}(M)$ and $\text{supp}(M)=\text{supp}(M^T) $.
	
	The uniqueness of $S$ is a consequence of Theorem \ref{RS}, and it is the limit of the Sinkhorn-Knopp-like algorithm by Theorem \ref{conv}. If $M$ is fully indecomposable its bipartite graph is connected \cite[Theorem 1.3.7]{Brualdi} and, thus, it is chainable \cite[Theorem 1.2]{Chainable} (see \cite{Chainable} or \cite{RoSc} for the definition of ``chainable''). Then, by \cite[Theorem 4]{RoSc}, $D_{M,\ell}$ and $D_{M,r}$ are also unique up to a nonnegative scalar multiple. \cblu{Finally, if, in this situation, $M=M^T$, then transposing both sides of $D_{M,\ell} M D_{M,r} \1_n = v$ and of $\1_n^T D_{M,\ell} M D_{M,r} = v^T$ implies $\1_n^T D_{M,r} M D_{M,\ell}  = v^T$ and $D_{M,r} M D_{M,\ell} \1_n = v$, which combined with the uniqueness of $D_{M,\ell}$ and $D_{M,r}$ up to an scalar multiple, implies that $D_{M,r} = \alpha D_{M,\ell}$ for some  $\alpha >0$, and $D = \sqrt{\alpha} D_{M,\ell}$ is the unique nonnegative diagonal matrix satisfying $S = DMD$.}
\end{proof}

If $M$ satisfies the conditions in Theorem \ref{th_squareRS}, the scaled matrix $S$ can be computed by using the Sinkhorn-Knopp-like algorithm in \cblu{Appendix $A$} with prescribed common vector $v$ for the row and column sums, i.e., with $r=c=v$.


In Section \ref{sec:regularized}, we will present new cost functions for our minimization problem \eqref{inf1} to make sure that it always has a unique and bounded solution. This new approach will be based on the results presented in this section combined with regularization techniques. In addition, this new approach will be applied to arbitrary pencils (square or nonsquare). First, we study in \cblu{Section \ref{sec:nonsquare} the unregularized} nonsquare case.

\section{Scaling nonsquare pencils and related problems} \label{sec:nonsquare}
In the square case, we scaled the pencil so that its row norms and column norms were equal as in \eqref{balanced}. However, this is no longer possible for $m\times n$ rectangular pencils since the numbers of rows and columns are different. But instead, one can try to balance the pencil $\la B-A$ by achieving the following equalities
\begin{equation}\label{columnrow_norms}
\begin{split}
& \|\text{col}_j (\widetilde{A})\|_2^2 + \|\text{col}_j (\widetilde{B})\|_2^2 = \gamma_\ell^2 \text{, for }j=1,\ldots, n, \text{ and} \\ & \|\text{row}_i (\widetilde{A})\|_2^2 + \|\text{row}_i (\widetilde{B})\|_2^2 = \gamma_r^2 \text{, for }i=1,\ldots, m,
\end{split}
\end{equation}
where $\widetilde{A}:=D_\ell A D_r$ and $\widetilde{B}:=D_\ell B D_r$ and $\|\la \widetilde{B}-\widetilde{A}\|_F^2=n\gamma_\ell^2=m\gamma_r^2$.
For the nonsquare case, we also define the nonnegative matrices
\begin{equation} \label{nonneg_nonsquare}
M :=  |A|^{\circ 2}+|B|^{\circ 2}, \quad \mathrm{and} \quad \widetilde{M} :=  |\widetilde{A}|^{\circ 2}+|\widetilde{B}|^{\circ 2}.
\end{equation}The scaling problem discussed in this section is a special case of the general scaling problem in Theorem \ref{RS}, where we choose $r=\gamma^2_r \1_m$ and
$c=\gamma^2_\ell \1_n$.

We now show that there is an optimization problem whose first order optimality conditions corresponds to the equalities in \eqref{columnrow_norms}.
\begin{theorem} \label{th:balanceAB}
	The following minimization problem over the set of positive diagonal matrices $D_\ell=\diag (d_{\ell_1}, \ldots ,d_{\ell_m})$ and $D_r=\diag (d_{r_1}, \ldots ,d_{r_n})$~:
	$$  \inf_{\det D_\ell^2=c_\ell,\det D_r^2 =c_r} (\| D_\ell A D_r \|_F^2 +\| D_\ell B D_r \|_F^2)
	$$
	has the balancing equations \eqref{columnrow_norms} as first order optimality conditions.
\end{theorem}

\begin{proof}
	If one makes the change of variables for the elements of $D_\ell$ and $D_r$ as follows $d^2_{\ell_i}=\exp(u_i)$,  $d^2_{r_j}=\exp(v_j)$,
	and introduce the notation $m_{ij}:=|a_{ij}|^2+|b_{ij}|^2$, then the above minimization is equivalent to a convex minimization problem with linear constraints~:
	\begin{equation}\label{constrained} \inf \sum_{i=1}^m \sum_{j=1}^n m_{ij} \exp(u_i+v_j), \quad \mathrm{subject \;\; to} \quad  \sum_{i=1}^m u_i = \ln c_\ell   , \quad \sum_{j=1}^n v_j= \ln c_r.
	\end{equation}
	The corresponding unconstrained problem with Lagrange multipliers $\Gamma_\ell$ and $\Gamma_r$, is
	$$ \inf \sum_{i=1}^m \sum_{j=1}^n m_{ij} \exp(u_i+v_j) + \Gamma_\ell(\ln c_\ell-  \sum_{i=1}^m u_i) + \Gamma_r(\ln c_r -  \sum_{j=1}^n v_j ).
	$$
	The first order conditions of optimality are the equality constraints of \eqref{constrained} and the equations
	\begin{equation}\label{eqconst}
	\sum_{j=1}^n d^2_{\ell_i}m_{ij}d^2_{r_j} = \Gamma_\ell, \quad  \sum_{i=1}^m d^2_{\ell_i}m_{ij}d^2_{r_j} = \Gamma_r,
	\end{equation}
	which express exactly that the row norms of $\widetilde M:= D_\ell^2MD_r^2$ are equal to each other and that its column norms are equal to each other.
	Since the Lagrange multipliers $\Gamma_\ell$ and $\Gamma_r$ are clearly nonnegative, we can can write them as  $\gamma^2_\ell:=\Gamma_\ell$ and $\gamma^2_r:=\Gamma_r$, which completes the proof.
\end{proof}

It is important to emphasize that unfortunately the optimization problem in Theorem \ref{th:balanceAB} does not always have a solution. This happens, for instance, if the corresponding matrix $M :=  |A|^{\circ 2}+|B|^{\circ 2}$ is the matrix appearing in Example \ref{ex:nonsquare}.

If there exists solution for the optimization problem in Theorem \ref{th:balanceAB}, it can be obtained by a sequence of alternating scalings $D_\ell^2$ and $D_r^2$ that make the rows of $D_\ell^2(MD_r^2)$ have equal sum $\gamma_r^2$,  and then the columns of $(D_\ell^2M)D_r^2$ have equal sum $\gamma_\ell^2$, while maintaining the  constraints
$\det D_\ell^2=c_\ell$, $\det D_r^2 =c_r$ \cblu{in the accumulated diagonal transformations, which determine the values of $\gamma_r^2$ and $\gamma_\ell^2$}. The cyclic alternation of row and column scalings, then amounts to coordinate descent applied to the minimization. This algorithm thus continues to decrease the cost function as long as the equalities \eqref{eqconst} are not met. \cblu{This is very similar to the Sinkhorn-Knopp-like Algorithm 1 applied to $M$ with $r = \gamma_r^2 \1_m$ and $c = \gamma_\ell^2 \1_n$. Since the exact values of $\gamma_r^2$ and $\gamma_\ell^2$ are of no interest, in practice one can simply apply Algorithm 1 to $M$ with $r = n \1_m$ and $c = m \1_n$.
Recall that,  according to Theorem \ref{conv}, this algorithm converges if and only if the corresponding scaling problem has solution.}

\begin{example}  \label{Ex1}
	Let us consider the pencil of a $5\times 6$ Kronecker block
	$$  \la B - A :=\left[\begin{array}{cccccc} \la & -1 \\  & \la & -1 \\ & & \la & -1 \\ & & & \la & -1 \\ & &  & & \la & -1
	\end{array}\right]
	$$
	then the scaled matrix $\widetilde M$ and the corresponding diagonal scaling matrices $D_\ell^2$ and $D_r^2$ look like
	\begin{equation} \label{eq.perfectM5x6}  \widetilde M :=\left[\begin{array}{cccccc} 5 & 1 \\  & 4 & 2 \\ & & 3  & 3 \\ & & & 2 & 4 \\ & &  & & 1 & 5
	\end{array}\right] ,  \quad \begin{array}{cc} D^2_\ell =\diag (1, 4, 6, 4, 1), & \gamma^2_\ell=5,  \\  \\ D_r^2=\diag (5, 1, 0.5, 0.5, 1, 5), & \gamma^2_r=6.\end{array}  \end{equation}
\end{example}

\section{The regularized scaling method for pencils}\label{sec:regularized}
The facts that for a nonsquare pencil the doubly stochastic scaling can not be applied anymore, that even for square pencils the corresponding matrix $M$ may not have total support and that the optimization problem in Theorem 4.1 does not always have solution can be by-passed by introducing a regularization term which will ensure an essentially unique bounded solution for $D_\ell$ and $D_r$. The cost of introducing such a term is that we will obtain a solution of an approximate problem. Nevertheless, with the new approach we can always assure that we will find such a solution.

Given two matrices $A,$ $B$ of size $m\times n,$ we consider the following constrained minimization problem over the set of \cblu{positive} diagonal matrices $D_\ell=\diag (d_{\ell_1}, \ldots ,d_{\ell_m})$ and $D_r=\diag (d_{r_1}, \ldots ,d_{r_n})$~:
\begin{equation}  \label{form1} \inf_{\det D_\ell^2\det D_r^2=c} 2(\| D_\ell A D_r \|_F^2 +\| D_\ell B D_r \|_F^2) + \alpha^2\left(\frac{1}{m^2}\|D_\ell \|_F^4+\frac{1}{n^2}\|D_r\|_F^4 \right),
\end{equation}
for some real number $c>0$ and a regularization parameter $\alpha$. If we denote again the matrix $M:=|A|^{\circ 2}+|B|^{\circ 2},$ then we can rewrite this as follows:
\begin{equation} \label{form2}  \inf_{\det D_\ell^2\det D_r^2 =c} \mathbf{1}_{m+n}^T
\left[\begin{array}{cc} \frac{\alpha^2}{m^2}D_\ell^2  \mathbf{1}_m \mathbf{1}_m^T D_\ell^2 & D_\ell^2 M D_r^2 \\
D_r^2 M^T D_\ell^2 &  \frac{\alpha^2}{n^2} D_r^2   \mathbf{1}_n \mathbf{1}_n^T D_r^2
\end{array}\right] \mathbf{1}_{m+n},
\end{equation}
which suggests that there may be a link to doubly stochastic scaling. Indeed, let us consider the two-sided scaling problem
$\widetilde M_\alpha := D_{\ell,r} M_\alpha D_{\ell,r}$, where
\begin{equation*} \label{total}
D_{\ell,r}:= \left[\begin{array}{cc} D_\ell  & 0 \\ 0&  D_r
\end{array}\right] ,
\end{equation*}
subject to $\det D_\ell^2\det D_r^2 = \det D^2_{\ell,r} = c,$ and
\begin{equation}\label{eq:Malpha}
M_\alpha^{\circ 2}=
\left[\begin{array}{cc} \frac{\alpha^2}{m^2} \mathbf{1}_m \mathbf{1}_m^T &  M  \\
M^T  &  \frac{\alpha^2}{n^2}  \mathbf{1}_n \mathbf{1}_n^T
\end{array}\right].
\end{equation} Notice that both diagonal blocks in $M_\alpha$ have Frobenius norm $\alpha$. We then prove in Theorem \ref{th_scaling} that the optimization problem  \eqref{form1} can be solved by the Sinkhorn--Knopp algorithm in a unique way. We will need the following auxiliary Lemma \ref{fullyindecomp} in our proof.

\begin{lemma}\label{fullyindecomp} Let $M_\alpha^{\circ 2}$ be the nonnegative matrix in \eqref{eq:Malpha} with $\alpha\neq0$. Then $M_\alpha^{\circ 2}$ has total support. Moreover, if $M \neq 0$ then $M_{\alpha}^{\circ 2}$ is fully indecomposable.
\end{lemma}
\begin{proof} See \cblu{Appendix B.}
\end{proof}

\begin{theorem}\label{th_scaling} Let $A$ and $B$ be $m\times n$ complex matrices and $\alpha,c > 0$ be real numbers. Let us consider the constrained minimization problem  \eqref{form1} over the set $\{(D_{\ell},D_{r}): D_\ell:=\diag (\d_{\ell_1},\ldots,\d_{\ell_m}), D_r:=\diag (\d_{r_1},\ldots,\d_{r_n}), \d_{\ell_i}>0, \d_{r_j}>0\}.$ Then the following statements hold:
	\begin{itemize}
		\item[a)] The optimization problem  \eqref{form1} is equivalent to the optimization problem  \eqref{form2}.
		\item[b)] The optimization problem  \eqref{form1} is equivalent to the optimization problem
		\begin{equation*}
		\inf_{\det D_\ell^2\det D_r^2 =c}  \left\|
		\left[\begin{array}{cc} D_\ell  & 0 \\ 0&  D_r
		\end{array}\right]
		M_\alpha  \left[\begin{array}{cc} D_\ell  & 0 \\ 0&  D_r
		\end{array}\right]\right\|_F^2 ,	\end{equation*}
		where $M_\alpha^{\circ 2}$ is given in \eqref{eq:Malpha}.
		\item[c)] There exists a unique solution $(\widetilde{D}_{\ell},\widetilde{D}_{r})$ of \eqref{form1}. Moreover, $(\widetilde{D}_{\ell},\widetilde{D}_{r})$ is bounded and makes the matrix
		$$\left[\begin{array}{cc} \widetilde{D}_\ell^2  & 0 \\ 0& \widetilde{ D}_r^2
		\end{array}\right]
		M_\alpha^{\circ 2}  \left[\begin{array}{cc} \widetilde{D}_\ell^2  & 0 \\ 0&  \widetilde{D}_r^2
		\end{array}\right]$$ a scalar multiple of a doubly stochastic matrix. Therefore, $(\widetilde{D}_{\ell},\widetilde{D}_{r})$ can be computed, \cblu{up to a scalar multiple}, by applying the algorithm in Appendix A to $M_\alpha^{\circ 2}$ \cblu{with $r=c=\1_{m+n}$}.
	\end{itemize}
\end{theorem}
\begin{proof} We have already seen statements $a)$ and $b)$ in this section because the optimization problem in $b)$ is just \eqref{form2}. Then we only need to prove $c).$ We make the change of variables $d^2_{\ell_i}=\exp(u_i)$ and $d^2_{r_j}=\exp(v_j)$ for the elements of $D_\ell$ and $D_r,$ respectively. Then the optimization problem \eqref{form1} is equivalent to  the optimization problem:
\begin{equation}\label{constrained2}
\begin{split}
	&\inf \;  2\sum_{i=1}^m \sum_{j=1}^n m_{ij} \exp(u_i+v_j)+\alpha^2\left(\frac{1}{m^2}\left(\sum_{i=1}^m \exp(u_i) \right)^2 + \frac{1}{n^2}\left(\sum_{j=1}^n \exp(v_j) \right)^2 \right),\\ & \mathrm{subject \;\; to} \quad  \sum_{i=1}^m u_i + \sum_{j=1}^n v_j = \ln c.
\end{split}	
\end{equation}
The corresponding unconstrained problem with Lagrange multiplier $\Gamma$ is:
\begin{equation}
\begin{split}
	\inf \; &  2\sum_{i=1}^m \sum_{j=1}^n m_{ij} \exp(u_i+v_j)+\alpha^2\left(\frac{1}{m^2}\left(\sum_{i=1}^m \exp(u_i) \right)^2 + \frac{1}{n^2}\left(\sum_{j=1}^n \exp(v_j) \right)^2 \right)\\ & + \Gamma \left(\ln c - \sum_{i=1}^m u_i - \sum_{j=1}^n v_j \right).
\end{split}	
\end{equation}
The first order conditions of optimality are the equality constraint of \eqref{constrained2} and the equations
\begin{equation*}
\frac{\alpha^2}{m^2}d^2_{\ell_i}\sum_{i=1}^m d^2_{\ell_i} + \sum_{j=1}^n d^2_{\ell_i}m_{ij}d^2_{r_j}	 = \dfrac{\Gamma}{2}, \quad \text{and} \quad 	\frac{\alpha^2}{n^2}d^2_{r_j}\sum_{j=1}^n d^2_{r_j} + \sum_{i=1}^m d^2_{\ell_i}m_{ij}d^2_{r_j}	 = \dfrac{\Gamma}{2},
\end{equation*}
for $i=1,\ldots,m$ and $j=1,\ldots,n$, respectively, which express that the row sum and the column sum of
$$\left[\begin{array}{cc} D_\ell^2  & 0 \\ 0&  D_r^2
	\end{array}\right]
	M_\alpha^{\circ 2}  \left[\begin{array}{cc} D_\ell^2  & 0 \\ 0&  D_r^2
	\end{array}\right]$$
are equal to $\dfrac{\Gamma}{2}.$ By Lemma \ref{fullyindecomp}, we know that $M_\alpha^{\circ 2}$ is fully indecomposable. Then, by the Sinkhorn--Knopp theorem, there exists a unique and bounded scaling $(E_{\ell},E_{r})$ that makes the matrix
$$\left[\begin{array}{cc} E_\ell^2  & 0 \\ 0& E_r^2
	\end{array}\right]
	M_\alpha^{\circ 2}  \left[\begin{array}{cc} E_\ell^2  & 0 \\ 0&  E_r^2
	\end{array}\right]$$
doubly stochastic. Assume that $\det E_\ell^2\det E_r^2=k.$ We define $\widetilde{D}_\ell :=\left(\frac{c}{k}\right)^{\frac{1}{2(m+n)}}E_{\ell}$ and $\widetilde{D}_r :=\left(\frac{c}{k}\right)^{\frac{1}{2(m+n)}}E_r.$ Then $\det \widetilde{D}_\ell^2\det \widetilde{D}_r^2=c$ and $(\widetilde{D}_{\ell},\widetilde{D}_{r})$ is the solution of \eqref{form1}. We can again redefine $\gamma^2:=\Gamma/2$ since this quantity is nonnegative.
\end{proof}

For completeness, we include the following result, which is a direct corollary of the proof of Theorem \ref{th_scaling}.

\begin{theorem}\label{th_scaling_rowcolumnsum} Let $A$ and $B$ be $m\times n$ complex matrices and $\alpha,c > 0$ be real numbers. Then the constrained minimization problem  \begin{equation*}   \inf_{\det D_\ell^2\det D_r^2=c} 2(\| D_\ell A D_r \|_F^2 +\| D_\ell B D_r \|_F^2) + \alpha^2\left(\frac{1}{m^2}\|D_\ell \|_F^4+\frac{1}{n^2}\|D_r\|_F^4 \right),
	\end{equation*} over the set $\{(D_{\ell},D_{r}): D_\ell:=\diag (\d_{\ell_1},\ldots,\d_{\ell_m}), D_r:=\diag (\d_{r_1},\ldots,\d_{r_n}), \d_{\ell_i}>0, \d_{r_j}>0\}$ has a unique and bounded solution. Moreover, it satisfies the equations:
	\begin{equation*}
	\begin{split}
	& \|\text{col}_j (\widetilde{A})\|_2^2 + \|\text{col}_j (\widetilde{B})\|_2^2+ \frac{\alpha^2}{n^2}\delta_{r_j}^2\|D_{r}\|_F^2= \gamma^2 \text{, for }j=1,\ldots, n, \text{ and} \\ & \|\text{row}_i (\widetilde{A})\|_2^2 + \|\text{row}_i (\widetilde{B})\|_2^2 + \frac{\alpha^2}{m^2}\delta_{\ell_i}^2\|D_{\ell}\|_F^2= \gamma^2 \text{, for }i=1,\ldots, m,
	\end{split}
	\end{equation*}
	for some nonzero scalar $\gamma $, where $\widetilde{A}:=D_\ell A D_r$ and $\widetilde{B}:=D_\ell B D_r$.
\end{theorem}

\begin{remark} \rm By Theorem \ref{th_scaling}, we know that the row sums and the column sums of the matrix
	$$\left[\begin{array}{cc} D_\ell^2  & 0 \\ 0&  D_r^2
	\end{array}\right]
	M_\alpha^{\circ 2}  \left[\begin{array}{cc} D_\ell^2  & 0 \\ 0&  D_r^2
	\end{array}\right]$$
	are equal to each other, where $(D_\ell, D_r)$ is the solution in Theorem \ref{th_scaling_rowcolumnsum}. The quantity of such row and column sums is the scalar $\gamma^2 $ appearing in Theorem \ref{th_scaling_rowcolumnsum}.
\end{remark}

\cblu{

In Example \ref{ex-f51}, we will illustrate the effect of choosing different values for the regularization parameter $\alpha$ in \eqref{eq:Malpha} in order to make the row and column sums of $D_\ell^2MD_r^2$ as equal as possible for a square matrix $M$ (corresponding to a pencil $\la B - A$) which does not have total support and, thus, cannot be scaled to a multiple of a doubly stochastic matrix. For measuring the quality of the obtained approximate scaling in this and other examples considered in this paper, we introduce the following definition.
\begin{definition} \label{def.quality-scaling}
Let $M \in \mathbb{R}^{m \times n}$ be a real nonnegative matrix, let $r(M) \in \mathbb{R}^{m \times 1}$ and $c(M) \in \mathbb{R}^{n\times 1}$ be, respectively, the vectors of row sums and column sums of $M$, denote by $r_i(M)$ and $c_i (M)$ the $i$-th entries of $r(M)$ and $c(M)$, and assume $r_i(M) > 0$ and $c_j (M)>0$ for all $i,j$. Then, the quality-factor of the homogeneous scaling of $M$ is defined as
\begin{equation} \label{eq.qs}
q_S (M) := \max \left\{ \frac{\max_i r_i (M)}{\min_i r_i (M)} \, , \, \frac{\max_i c_i (M)}{\min_i c_i (M)}  \right\} \, .
\end{equation}
\end{definition}
Observe that $q_S(M) = 1$ if and only if the row sums of $M$ are all equal and the column sums of $M$ are all equal. The closer to $1$ the factor $q_S(M)$ is, the better balanced the matrix $M$ is.

\begin{example} \label{ex-f51} We consider the square pencil $\la B_{1} - A_{1}$ in Example \ref{examples}. The associated matrix
	\begin{equation} \label{eq.matex51}
	M_1 :=\left[\begin{array}{ccc} 1 & 1 & 0 \\ 1 & 0 & 0 \\ 0 & 0 & 1 \end{array}\right]
	\end{equation}
	has no total support and, thus, the Sinkhorn-Knopp algorithm does not converge. More precisely, the algorithm in Appendix A applied to $M_1$ with $r = c = \1_3$ and {\tt tol}$= 10^{-3}$ does not converge after $1000$ steps. In contrast, the same algorithm applied to the matrix $M_\alpha^{\circ 2}$ in \eqref{eq:Malpha} with $r=c=\1_6$ and {\tt tol}$= 10^{-3}$ for $\alpha = 1,0.5,0.1$ converges and produces scaled matrices $\widetilde M_\alpha = D_\ell^2MD_r^2$ which are approximately doubly stochastic up to a scalar multiple. The results are shown in Table \ref{table.ex51}, where the last column shows the 2-norm condition numbers of $D^2_\ell \approx D_r^2 $ and {\rm steps} denotes the number of steps until convergence, with each step comprising one right and one left diagonal scaling.
\begin{table}[h!]
\cblu{ \begin{center}
    \caption{Results of the regularization applied to the matrix in \eqref{eq.matex51} with {\tt tol}$= 10^{-3}$. The quality factors $q_S$ should be compared with $q_S(M_1) = 2$}
    \label{table.ex51}
    \begin{tabular}{c|c|c|c|c}
    $\alpha$ & steps & $q_S$ & $\diag (D^2_\ell) \approx \diag(D_r^2) $
    &  $\kappa (D^2_\ell) \approx \kappa (D_r^2) $
    \\ \hline
    1  &   11  &   1.38 &  0.485  , 1.29   ,   0.864  & 2.66 \\
    0.5  &  24  & 1.19 & 0.395    ,   2.05  ,   0.952 & 5.19 \\
    0.1  &  124 & 1.04 &  0.187  ,     5.15   ,   0.970 & 27.5
    \end{tabular}
\end{center} }
\end{table}

\noindent Choosing a smaller $\alpha$ yields a better equilibration for the row and column sums as measured by the quality-factor $q_S$ (to be compared with $q_S (M_1) =2$ for the original matrix), but at the cost of a worse conditioning of the scaling matrices $D^2_\ell, D_r^2$  and of a slower convergence. The latter is to be expected since for $\alpha=0$ the scaling to a multiple of a double stochastic matrix does not exist for $M_1$.

Finally, we show the results obtained when the algorithm in Appendix A is applied directly to $M_1$ with $r = c = \1_3$, i.e., without any regularization, but with the very relaxed stopping criterion {\tt tol}$= 1$. In this case the algorithm converges in only $3$ steps and the results are shown in Table \ref{table.ex51-bis}, where $\alpha = 0$ indicates that the problem has not been regularized (though the matrix $M_\alpha^{\circ 2}$ is not used at all). We will use this convention in other numerical examples and tests.
\begin{table}[h!]
\cblu{ \begin{center}
    \caption{Results of the unregularized Sinkhorn-Knopp algorithm applied to the matrix in \eqref{eq.matex51} with {\tt tol}$= 1$. The quality factor $q_S$ should be compared with $q_S(M_1) = 2$}
    \label{table.ex51-bis}
    \begin{tabular}{c|c|c|c|c|c|c}
    $\alpha$ & steps & $q_S$ & $\diag (D^2_\ell)$ & $\diag(D_r^2)$
    &  $\kappa (D^2_\ell)$ & $\kappa (D_r^2) $
    \\ \hline
    0  &  3 & 1.33 &  0.350, 2.45 , 0.765 & 0.408, 2.45, 1.31 & 7 & 6
    \end{tabular}
\end{center} }
\end{table}
The motivation for computing this rough {\tt tol}$= 1$ approximate solution will be clear in Section \ref{sec:numerics} and is related to the fact, previously commented, that for the purpose of improving the accuracy of the eigenvalues of $\la B -A$ computed in floating point arithmetic it is essential that the entries of the diagonal scaling matrices $D_\ell$ and $D_r$ are integer powers of $2$. This implies that it makes no sense to compute very precise scaling matrices $D_\ell$ and $D_r$, since their entries will be later rounded to their nearest integer powers of $2$ and, thus, a relaxed stopping criterion can be used. We remark here three facts that will be further discussed in Section \ref{sec:numerics}: {\tt tol}$= 1$ very often has a regularization effect, speeds up considerably the convergence and yields a reasonably ``well balanced'' matrix.
\end{example}

\begin{remark} \label{rem.regular} The choice of the regularization parameter $\alpha$ has to be guided by the equilibrium one wants to achieve between the ``quality'' of the balancing, the boundedness/conditioning of the diagonal scaling matrices and the speed of convergence. This depends heavily on the applied problem the user wants to solve. For the problem of improving the accuracy of computed eigenvalues, we do not need to consider a very small value of $\alpha$ since, in practice, it is enough to get a reasonably ``well balanced'' matrix $M$, because the entries of the diagonal scaling matrices have to be later rounded to their nearest integer powers of two. Moreover, as we will see in Section \ref{sec:numerics}, the use of the relaxed stopping criterion {\tt tol}$= 1$ makes it often unnecessary the use of the regularization. This can happen even in cases where the use of the regularization is mandatory from a theoretical point of view, since there is no exact solution of the scaling problem. The use of {\tt tol}$= 1$ prevents, in any case, to obtain very ``well-balanced matrices''. Thus, for the eigenvalue problem, we recommend to start always by using the un-regularized method and if it does not converge in a small number of steps (say $n/10$ for large $n$) to change to the regularized method with a value of $\alpha \lesssim 0.5 \max_{ij} \sqrt{M_{ij}}$. In contrast, in other type of problems where it is important to get always a very ``well-balanced matrix'' and a relaxed stopping criterion is not adequate or neccessary, a recommendable option might be to always use the regularization with a small value of $\alpha$, especially when $M$ is sparse, since it guarantees the existence of a solution. This will increase the complexity of the Sinkhorn-Knopp algorithm each step by a factor 4 since the matrix sizes are doubled. In difficult cases, this might be very slow and, thus, the regularized problem and the Sinkhorn-Knopp algorithm should be combined with faster algorithms (see \cite{Idel,CuturiPeyre} for the state-of-the art).
\end{remark}
}

\begin{remark} One could also have considered for the regularization the cost function
	$$ \inf_{\det D_\ell^2\det D_r^2=c} 2(\| D_\ell A D_r \|_F^2 +\| D_\ell B D_r \|_F^2) + \alpha^2\left(\|D_\ell^2 \|_F^2+\|D_r^2\|_F^2 \right),$$
	which would correspond to the matrix
	$$ M_\alpha^{\circ 2} := \left[\begin{array}{cc}  \alpha^2I_m &  M  \\
	M^T  &   \alpha^2 I_n \end{array}\right].$$
	This matrix has total support for $\alpha>0$. However, it is not necessarily fully indecomposable (assume for instance that $M$ has a zero row or column)  and, therefore, we can not guarantee the essential uniqueness of the scaling matrices $D_\ell$ and $D_r$.
\end{remark}

\subsection{The regularized method with prescribed nonhomogeneous common vector for the row and column sums}\label{sec:regularized_rowandcolum}

In the nonsquare case, we know from the discussions of Section \ref{sec:nonsquare} that making the column and row sums of $\widetilde M=D_\ell^2MD_r^2$ become equal can not be achieved exactly, where $M$ is the matrix in \eqref{nonneg_nonsquare}. In this case, we can use the regularized method in Theorem \ref{th_scaling}$-c)$ in order to obtain a scaling that balances $\widetilde M$ approximately. We have used this approach on many problems and have obtained pretty satisfactory results. However, since by using this method we always obtain a scalar multiple of a doubly stochastic matrix as solution for $M_\alpha^{\circ 2}$, this method considers in some sense the rows and columns of $M$ in the same way, which is not natural in the rectangular case. Thus, one possible strategy for improving this approach is not to request that $M_\alpha^{\circ 2}$ is scaled to be a scalar multiple of a doubly stochastic matrix but to impose a modified scaling with prescribed common vector
\begin{equation} \label{eq.nonhomogv}
v:=\left[\begin{array}{c} n\mathbf{1}_m \\ m\mathbf{1}_n \end{array}\right]
\end{equation}
for the row and column sums. The new regularized method is then described by :
\begin{equation}\label{eq:rowsums}
\left[\begin{array}{cc} D_\ell^2  & 0 \\ 0&  D_r^2
\end{array}\right]
M_\alpha^{\circ 2}  \left[\begin{array}{cc} D_\ell^2  & 0 \\ 0&  D_r^2
\end{array}\right] \left[\begin{array}{c} \mathbf{1}_m \\ \mathbf{1}_n \end{array}\right] =v
\end{equation}
and
\begin{equation}\label{eq:columsums}
\left[\begin{array}{cc} \mathbf{1}_m^T & \mathbf{1}_n^T \end{array}\right] \left[\begin{array}{cc} D_\ell^2  & 0 \\ 0&  D_r^2
\end{array}\right]
M_\alpha^{\circ 2}  \left[\begin{array}{cc} D_\ell^2  & 0 \\ 0 &  D_r^2
\end{array}\right] = v^{T}.
\end{equation}

This is a problem that falls into the category of scalings considered in Theorem \ref{RS}. In addition, notice that the matrix $M_\alpha^{\circ 2}$ satisfies the hypotheses in Theorem \ref{th_squareRS} if $\alpha\neq 0$, i.e., $\text{supp}(M_\alpha^{\circ 2})=\text{supp}((M_\alpha^{\circ 2 })^{T}) $ and $(i,i)\in\text{supp}(M_\alpha^{\circ 2})$ for all $i=1,\cdots,n+m$. Then, by considering $\alpha\neq 0$, we know by Theorem \ref{th_squareRS} that there always exists a solution for this modified scaling problem with prescribed common vector for the row and column sums. Moreover, since $M_\alpha^{\circ 2}$ is fully indecomposable \cblu{when $M \ne 0$, according to Lemma \ref{fullyindecomp}, and is symmetric, there exists a unique and bounded diagonal scaling matrix $\diag(D_\ell^2 ,D_r^2)$ solving the problem \eqref{eq:rowsums}-\eqref{eq:columsums}, according again to Theorem \ref{th_squareRS}}. It can also be computed by using the Sinkhorn-Knopp-like algorithm \cblu{given in Appendix A with $r = c = v$,  as it converges to the unique solution by Theorem \ref{th_squareRS}. In our numerical experience, this approach very often improves, for rectangular matrices $M$, the results with respect to the approach in Theorem \ref{th_scaling}$-c)$ (corresponding to apply to $M_\alpha^{\circ 2}$ the algorithm in Appendix A with $r = c = \1_{m+n}$) in terms of the number of steps until convergence and the quality of the scaling of the obtained matrix.}


Notice that, when $\alpha=0$, \cblu{the scaling problem \eqref{eq:rowsums}-\eqref{eq:columsums}} reduces to the problem discussed in Section \ref{sec:nonsquare}. Then, for very small $\alpha$, the regularized scaling with prescribed row and column sums $v$ tends to the scaling problem explained in Section \ref{sec:nonsquare}, which does not always have a solution.

In the following example, we illustrate the effect of choosing different values of $\alpha$ and the row and column sum conditions \eqref{eq:rowsums} and \eqref{eq:columsums}.

\begin{example} \label{Ex2} \cblu{We remark that, for this example, the algorithm described in Section \ref{sec:nonsquare} converges. More precisely, the algorithm in Appendix A applied to the matrix $M = |A|^{\circ 2}+|B|^{\circ 2}$ with $r = 6\cdot \1_5$ and $c = 5 \cdot \1_6$ converges. Thus, there is no need to use the regularized method. Nevertheless, we use the regularized method developed in this section with two purposes: (1) for comparing the approximate regularized solution and the exact solution of the optimization problem in Theorem \ref{th:balanceAB} and (2) for illustrating the effect of choosing different values of $\alpha$.
We consider again the nonsquare pencil $\lambda B - A$ in Example \ref{Ex1} but now with a preliminary diagonal scaling
	$\lambda \hat B - \hat A:=\hat D_\ell (\lambda B-A)\hat D_r$ on the left and the right with condition numbers  $\kappa(\hat D_\ell) = 12.3$ and $\kappa(\hat D_r)= 2409.1$. The resulting matrix $M:=\hat A^{\circ 2}+\hat B^{\circ 2} $ to be scaled is
{\small \begin{equation} \label{eq.matorigex52}
M = \left[ \begin{array}{cccccc}
 8.983e-06  & 1.145e-09   &         0   &         0   &         0   &         0 \\
            0 &  1.231e-08 &  6.801e-02  &          0        &    0      &      0 \\
            0 &           0 &  4.734e-02  & 1.228e-02      &      0    &        0 \\
            0 &           0 &           0  & 1.977e-03 &  5.170e-04    &        0 \\
            0 &           0 &           0  &    0 &  6.464e-02  & 1

\end{array}\right]
\end{equation} }
which we normalized to have its largest element equal to $1$. This is a severely unbalanced matrix with quality-factor $q_S (M) = 7.43 \cdot 10^7$, as defined in \eqref{eq.qs}, which combined with the sparsity of the matrix, makes it a difficult problem for the Sinkhorn-Knopp-like algorithm. When applying to $M$ the algorithm in Appendix A with $r = 6\cdot \1_5$, $c = 5 \cdot \1_6$ and {\tt tol}$= 10^{-3}$, we obtained (with three digits of accuracy)
the same result as in Example \ref{Ex1}, i.e., the matrix in \eqref{eq.perfectM5x6}. This indicates that the direct scaling method can compensate for a bad initial scaling. The other results of this unregularized method are displayed in the first line of Table \ref{table.ex53}.

We now apply the regularized method with the matrix $M_\alpha^{\circ 2}$ and prescribed common vector $v:=[6,6,6,6,6,5,5,5,5,5,5]^{T}$ for the row and column sums, i.e., the algorithm in Appendix A applied to $M_\alpha^{\circ 2}$ with $r = c = v$ and {\tt tol}$= 10^{-3}$, for three different values of $\alpha$. The results are shown in Table \ref{table.ex53}. These results show that decreasing $\alpha$ in the regularized method improves the quality of the scaling, but makes the diagonal scaling matrices worse conditioned and the convergence slower. Also one can see that the regularization yields considerable improvements of the scaling with respect to the original matrix $M$ with not too small $\alpha$ and with a comparable number of steps to the regularized method (see, for instance, the results for $\alpha = 10^{-3}$). However, the convergence of the regularized method to the unregularized one when $\alpha \rightarrow 0$ is slow. In this example $\alpha = 10^{-5}$ and $steps = 1147$ are needed to get $q_S = 1.004$ with {\tt tol}$= 10^{-3}$.

\begin{table}[h!]
\cblu{ \begin{center}
    \caption{Results of the unregularized and the regularized methods applied to the matrix in \eqref{eq.matorigex52} with $r = 6\cdot \1_5$ and $c = 5 \cdot \1_6$ for the unregularized method, $v$ as in \eqref{eq.nonhomogv} for the regularized method and {\tt tol}$= 10^{-3}$ in all cases. The quality factors $q_S$ should be compared with $q_S (M) = 7.43 \cdot 10^7$}
    \label{table.ex53}
    \begin{tabular}{c|c|c|c|c}
    $\alpha$ & steps & $q_S$ & $\kappa (D^2_\ell)$ & $\kappa (D_r^2) $
    \\ \hline
     0      &     94  &     1   &     499.3  & 1.1066e+07 \\
     $10^{-2}$   &        53  &     844   &    1204.6  &     7330.3 \\
     $10^{-3}$   &        99  &     8.04   &    766.4  & 8.5339e+05 \\
     $10^{-4}$   &       154  &     1.20   &    501.4  & 7.4396e+06  \\
    \end{tabular}
\end{center} }
\end{table}

Finally, as in Example \ref{ex-f51} and based on the same motivations explained there, we show in Table \ref{table.ex53-bis} the results of applying directly to $M$ in \eqref{eq.matorigex52}, the algorithm in Appendix A with $r = 6\cdot \1_5$, $c = 5 \cdot \1_6$ and the relaxed stopping criterion {\tt tol}$= 1$. The results are extremely good in terms of the speed of convergence and the improvement of the quality of the scaling $q_S$.
\begin{table}[h!]
\cblu{ \begin{center}
    \caption{Results of the unregularized Sinkhorn-Knopp-like algorithm applied to the matrix in \eqref{eq.matorigex52} with $r = 6\cdot \1_5$, $c = 5 \cdot \1_6$ and {\tt tol}$= 1$. The quality factor $q_S$ should be compared with $q_S (M) = 7.43 \cdot 10^7$}
    \label{table.ex53-bis}
    \begin{tabular}{c|c|c|c|c}
    $\alpha$ & steps & $q_S$ & $\kappa (D^2_\ell)$ & $\kappa (D_r^2) $
    \\ \hline
     0      &     4  &     1.59  &     119.97  &  2.8292e+07 \\
    \end{tabular}
\end{center} }
\end{table}

}
\end{example}

\cblu{As commented in Sections \ref{sec:scaling} and \ref{sec:regular}, in the rectangular case, simple necessary and sufficient conditions on the zero pattern of $M$ for the scaling technique in Section \ref{sec:nonsquare} (i.e., the algorithm in Appendix $A$ applied to $M$ with $r = n \1_m$ and $c = m \1_n$) to converge are not known (see \cite{Idel,CuturiPeyre} for the state of the art)}. In contrast, the regularized method with the matrix
$M_\alpha^{\circ 2}$ and prescribed common vector $v$ in \eqref{eq.nonhomogv} always has a solution for rectangular pencils, and the previous example, as well as many others, shows that it produces satisfactory results, \cblu{even when the unregularized problem has solution}. Therefore, using this new regularized method is \cblu{always an available option for scaling a rectangular $M$, regardless of whether the optimization problem} in Theorem \ref{th:balanceAB} has a solution or not.

In Example \ref{Ex2}, we knew that the corresponding matrix $M$ can be scaled with prescribed vectors $r:=[6,6,6,6,6]^{T}$, for the row sums, and  $c:=[5,5,5,5,5,5]^{T}$, for the column sums. We now consider the matrix $M$ in Example \ref{ex:nonsquare} that can not be scaled \cblu{to have equal row sums and equal column sums}, but we use the regularized method with prescribed common vector \eqref{eq.nonhomogv} for the row and column sums to obtain an approximate scaling.

\cblu{
\begin{example} We consider the nonsquare matrix
	\begin{equation} \label{eq.matex54}
	M:=\left[\begin{array}{ccc} 1 & 1 & 1 \\ 0 & 0 & 1\end{array}\right]
	\end{equation}
	in Example \ref{ex:nonsquare}, that can not be scaled with prescribed vectors $r:=[3,3]^{T}$, for the row sums, and  $c:=[2,2, 2]^{T}$, for the column sums. Therefore, the algorithm in Section \ref{sec:nonsquare}, i.e., the algorithm in Appendix A with this $r$ and $c$, does not converge for this matrix, neither with a stringent stopping criterion {\tt tol}$= 10^{-3}$ nor with the relaxed one {\tt tol}$= 1$ (which shows that {\tt tol}$= 1$ does not always yield convergence). More precisely, we have run this algorithm until $10^4$ steps and it gets stuck, alternating periodically between the following two matrices
\[
M_{\infty,1} = \left[ \begin{array}{ccc}
       1.5 & 1.5 & 4.9407e-324 \\
       0 & 0 & 3
     \end{array} \right] \quad \mbox{and} \quad
M_{\infty,2} = \left[ \begin{array}{ccc}
       2 & 2 & 4.9407e-324 \\
       0 & 0 & 2
     \end{array} \right]    .
\]
Observe that the quality-factors for the homogeneous scalings of the three matrices above are $q_S (M) = 3$ and $q_S (M_{\infty,1}) =q_S (M_{\infty,2}) = 2$, which means that although the un-regularized method does not converge, it has progressed towards a better scaling. Then, we use the regularized approach with different values of $\alpha$ and prescribed common vector $v:=[3,3,2,2,2]^{T}$ for the row and column sums of $M_{\alpha}^{\circ 2}$, i.e., the algorithm in Appendix A applied to $M_{\alpha}^{\circ 2}$ with $r = c = v$ and {\tt tol}$= 10^{-3}$. The results are shown in Table \ref{table.ex54}, where we observe that the regularization yields, even for rather large values of $\alpha$, a significant improvement in the quality of the scaling with a moderate number of steps and well-conditioned $D_\ell$ and $D_r$. In our experiment, $q_S$ reaches quickly a limit value of $1.5$ as $\alpha \rightarrow 0$ with the following corresponding limiting scaled matrix for $\alpha = 10^{-10}$:
\[
M_{\alpha \rightarrow 0} = \left[ \begin{array}{ccc}
 1.5    &      1.5  & 1.0301e-20 \\
            0    &        0     &       2
     \end{array} \right].
\]

\begin{table}[h!]
\cblu{ \begin{center}
    \caption{Results of the regularized method applied to the matrix in \eqref{eq.matex54} with $v:=[3,3,2,2,2]^{T}$ and {\tt tol}$= 10^{-3}$ in all cases. The quality factors $q_S$ should be compared with $q_S (M) = 3$}
    \label{table.ex54}
    \begin{tabular}{c|c|c|c|c}
    $\alpha$ & steps & $q_S$ & $\kappa (D^2_\ell)$ & $\kappa (D_r^2) $
    \\ \hline
    0.5   &  14    &   1.6441  &      10.39    &   8.0413 \\
    $10^{-1}$  & 20 & 1.5073 & 198.27 & 148.92 \\
    $10^{-2}$  & 29 & 1.5001 & 19422 & 14566  \\
    $10^{-4}$  & 45 & 1.5 & 1.9416e+08 & 1.4562e+08  \\
    $10^{-10}$ & 93 & 1.5 & 1.9416e+20 & 1.4562e+20
    \end{tabular}
    \end{center} }
\end{table}
\end{example}
}

\cblu{We end this section by looking} at the effect of the two sided scaling on the sensitivity of the underlying eigenvalue problem.
In the case of regular pencils, we argued \cite{LemVD} \cblu{(see also the discussion in  Section \ref{sec:regular})} that the minimization problem
$$  \inf_{\det T_\ell \det T_r =1} \| T_\ell (\la B - A) T_r \|_F^2 ,$$
over the arbitrary nonsingular matrix pairs $(T_\ell,T_r)$, yielded nearly optimal sensitivity for the generalized eigenvalues of the pencil.
But since the eigenvalue problem for a singular pencil is known to be ill-conditioned, this may not make sense anymore. Nevertheless, if we constrain the transformations to be bounded, then the Kronecker structure can not change anymore, and it \cblu{then makes sense} to talk about the sensitivity of the eigenvalues again. In the numerical examples we show that the scaling also improves the sensitivity of the eigenvalues of the regular part of a singular pencil.

\section{Numerical examples} \label{sec:numerics}
\cblu{In this section, we verify in many numerical tests that the scaling procedures described in Sections \ref{sec:regular}, \ref{sec:nonsquare} and \ref{sec:regularized} indeed improve the accuracy of computed eigenvalues of arbitrary pencils with a much smaller cost than computing the eigenvalues by the $QZ$ or staircase algorithms \cite{QZ,Van79}. All the numerical tests in this paper were performed in MATLAB R2019a. In Subsection \ref{sub:1}, we focus on the computational cost of the scaling procedures, which is much smaller than the cost of computing the eigenvalues as a consequence of the use of the relaxed stopping criterion {\tt tol}$=1$ in the algorithm in Appendix A. In Subsection \ref{sub:2}, we compare the accuracy of the computed eigenvalues of regular pencils without scaling and after the scaling described in Section \ref{sec:regular}. Moreover, we also compare the results with those corresponding to the scaling method of Ward \cite{Ward}, which is the only method currently implemented in LAPACK for scaling regular pencils\footnote{\cblu{Neither MATLAB nor LAPACK \cite{lapack} include built-in functions or routines for computing eigenvalues of singular pencils.}}. This comparison was already performed in \cite{LemVD} but only for regular pencils of dimension $10 \times 10$. Our experiments confirm that the method described in Section \ref{sec:regular}, i.e., that in \cite{LemVD}, outperforms Ward's method, which has a very poor behavior for certain pencils. In Subsection \ref{sub:3}, we perform similar tests on square singular pencils applying either the un-regularized scaling in Section \ref{sec:regular} or, if necessary, the regularized one in Section \ref{sec:regularized} and extract similar conclusions. Finally, in Subsection \ref{sub:secrectangularnum}, we perform tests on rectangular pencils applying either the un-regularized scaling in Section \ref{sec:nonsquare} or, if necessary, the regularized one in Subsection \ref{sec:regularized_rowandcolum}, which improve significantly the accuracy of the computed eigenvalues.}



\subsection{The stopping criterion {\tt tol}$=1$, computational cost and regularization} \label{sub:1}
\cblu{Given a complex $m\times n$ pencil $\la B - A$, all the scaling procedures described in this paper start by constructing the nonnegative matrix $M :=  |A|^{\circ 2}+|B|^{\circ 2}$. Then, the unregularized methods in Sections \ref{sec:regular} and \ref{sec:nonsquare} apply the algorithm in Appendix A to $M$ with $r = n \1_m$ and $c = m \1_n$, which in the square case means $r = c = n \1_n$. On the other hand, the regularized methods in Section \ref{sec:regularized} apply the algorithm in Appendix A to the nonnegative matrix $M_\alpha^{\circ 2}$ in \eqref{eq:Malpha} with $r = c = (2n) \1_{2n}$, when $m=n$, or $r = c = v$ in the rectangular case, where $v$ is the vector in \eqref{eq.nonhomogv}. In both, the unregularized and the regularized methods, one obtains a scaled matrix $\widetilde{M}= D_\ell^2 M  D_r^2$, together with the diagonal matrices $D_\ell^2$, $D_r^2$. Then, the scaling process of the pencil finishes in exact arithmetic by computing $D_\ell$, $D_r$, $\widetilde{A}= D_\ell A D_r$ and $\widetilde{B}= D_\ell B D_r$, with the aim of computing the eigenvalues of $\la \widetilde{B} - \widetilde{A}$ via some numerical algorithm. However, in real practice this must be applied in a computer and, then, there are rounding errors in the computation of $\widetilde{A}= D_\ell A D_r$ and $\widetilde{B}= D_\ell B D_r$. This implies that the pencils $\la B - A$ and $\la \widetilde{B} - \widetilde{A}$ are not exactly strictly equivalent to each other and, in the case $D_\ell$ and $D_r$ are ill conditioned as often happens in practice, their eigenvalues may be very different to each other and the process would not be useful for improving the accuracy of computed eigenvalues. In the spirit of the classical reference \cite{Par} (see also \cite{LemVD,Ward}), we can circumvent this difficulty if once $D_\ell$ and $D_r$ have been computed, we replace their diagonal entries by their nearest integer powers of $2$ to get new $D_\ell$ and $D_r$. With these new approximate diagonal scalings, $\widetilde{A}= D_\ell A D_r$ and $\widetilde{B}= D_\ell B D_r$ are computed exactly in floating point arithmetic and $\la B - A$ and $\la \widetilde{B} - \widetilde{A}$ have exactly the same eigenvalues. Of course, in this way, we do not obtain the same scaled pencil as in exact arithmetic, but it is expected that the obtained one is good enough for improving the accuracy of the computed eigenvalues.

The discussion above indicates that for eigenvalue computations, it is not needed to apply the algorithm in Appendix A to either $M$ or $M_\alpha^{\circ 2}$ with a stringent stopping criterion, because we will replace anyway the entries of $D_\ell$ and $D_r$ by their nearest integer powers of $2$. The stopping criterion of the algorithm in Appendix A applied to $M$ used for {\em the updating scaling} $D_{\ell,up}$ and $D_{r,up}$ in the iterative procedure is
$$
\max \left\{ 1- \frac{1}{\kappa(D^2_{\ell,up})} , 1- \frac{1}{\kappa(D^2_{r,up})} \right\} < \frac{\mathtt{tol}}{2}
$$
in terms of the spectral condition numbers of $D^2_{\ell,up}$ and  $D^2_{r,up}$. This is equivalent to
$$
\max \left\{ \kappa(D^2_{\ell,up}) \, , \, \kappa(D^2_{r,up}) \right\} < 1 + \frac{\mathtt{tol}}{2- \mathtt{tol}} \, .
$$
Thus,  {\tt tol}$=1$ implies that the algorithm stops when both $D_{\ell,up}$ and $D_{r,up}$ have a condition number smaller than $\sqrt{2}$. Since we are
approximating the final scaling matrices to integer powers of 2, this is a safe stopping criterion for practical purposes.

We will use {\tt tol}$=1$ in all the experiments in Subsections \ref{sub:2}, \ref{sub:3} and \ref{sub:secrectangularnum}. In the rest of this subsection, we will present some numerical tests that illustrate the impact of {\tt tol}$=1$ on the reduction of the number of steps that the algorithm in Appendix A needs for convergence and on the regularization of the problem. In all the tables for the experiments in this section ``steps'' denotes the number of steps until convergence, where one step includes one right and one left diagonal scaling. Moreover, $q_S (M_{orig})$ denotes the quality-factor defined in \eqref{eq.qs} for the original matrix $M =  |A|^{\circ 2}+|B|^{\circ 2}$ and $q_S (M_{scal})$ the one of the scaled matrix\footnote{\cblu{We emphasize that in all the experiments in Section \ref{sec:numerics}, the matrix $\widetilde{M}$ is computed as $\widetilde{M} = |D_\ell A D_r|^{\circ 2} + |D_\ell B D_r|^{\circ 2}$, where the diagonal matrices $D_\ell$ and $D_r$ are the ones whose diagonal entries are integer powers of $2$.}} $\widetilde{M}= D_\ell^2 M  D_r^2$. The ideal goal of all our scalings procedures is to make the row sums of $\widetilde{M}$ as equal as possible and its column sums as equal as possible as well, i.e, to get $q_S (M_{scal}) \approx 1$. As discussed in previous sections, we know that this is not always possible in exact arithmetic. In addition, even when it is possible in exact arithmetic, the use of entries that are integer powers of $2$ in $D_\ell$ and $D_r$ prevents to get such a goal. Thus, the practical goal is to get that $q_S (M_{scal})$ is much closer to $1$ than $q_S (M_{orig})$.

In our first test, we chose pencils of dimension $n\times n$ with $n=400, 800, 1200, 1600,$ $2000$ and with elements that were generated using MATLAB's {\tt randn} function elevated to power 20, yielding matrices $M$ with row and column sums strongly unbalanced. For each size $n$, we ran the algorithm in Appendix A with $r = c = n \1_n$ on ten random pencils and averaged the different tested magnitudes, both with {\tt tol}$= 1$ and {\tt tol}$= 10^{-3}$ and, in both cases, approximating $D_\ell$ and $D_r$ by their nearest integer powers of $2$. The results are shown in Table \ref{tab:table61-fro}. We emphasize that in this test, regularization is not needed because the random generation used for $A$ and $B$ imply that the entries of $M$ are almost always different from zero and, thus, $M$ has total support. Observe, that {\tt tol}$= 1$ yields a much faster convergence and similar values of $q_S (M_{scal})$ than {\tt tol}$= 10^{-3}$, which is very slow on this highly unbalanced matrices. Moreover, the number of required iteration steps does not grow with the dimension of the pencils. Since each step of the scaling procedure costs $O(n^2)$ flops, while the cost of computing the eigenvalues of an $n \times n$ pencil with the $QZ$ algorithm is $30 n^3$ flops  \cite[Section 7.7]{GoVL}, we conclude that for the matrices in this test the computational cost of the scaling procedure with {\tt tol}$= 1$ is much smaller than the cost of computing the eigenvalues.

}

\begin{table}[h!]
	\begin{center}
\cblu{
		\caption{Numerical test illustrating that the use of {\tt tol}$= 1$ decreases very much the number of steps without affecting to the quality of the scaling of the achieved scaled matrix $\widetilde{M}$ nor to the condition numbers of $D_\ell$ and $D_r$. The algorithm in Appendix A with $r = c = n \1_n$ has been applied to the matrices $M$ of exactly the same set of $n\times n$ pencils generated in {\rm MATLAB} as {\tt A=randn(n,n).}\!\!\^ {\tt (20)} and {\tt B=randn(n,n).}\!\!\^{\tt (20)}, one time with {\tt tol}$= 1$ and another time with {\tt tol}$= 10^{-3}$. No regularization is used, which is indicated with  $\alpha =0$}
		\label{tab:table61-fro}
	{\small	
\begin{tabular}{c|c||c|c|c|c||} 
  & & \multicolumn{4}{c||}{{\tt tol}$= 1$ and $\alpha =0$}  \\ \hline		
	$n$ & $q_S(M_{orig})$ & $q_S (M_{scal})$ & $\kappa (D_\ell)$ & $\kappa (D_r)$ & steps  \\	\hline
 400 & 1.94e+10 & 1.24e+01 & 2.76e+03 & 4.30e+04 & 9.8  \\
 800 & 4.90e+09 & 1.37e+01 & 2.46e+03 & 2.54e+04 & 10  \\
 1200 & 1.12e+10 & 1.35e+01 & 2.97e+03 & 1.35e+04 & 10.9  \\
 1600 & 2.79e+09 & 1.37e+01 & 2.56e+03 & 1.23e+04 & 10.7   \\
 2000 & 4.07e+09 & 1.42e+01 & 2.00e+03 & 1.37e+04 & 10.8
\end{tabular}

\vspace*{0.3cm}

\begin{tabular}{c|c||c|c|c|c||} 
  & & \multicolumn{4}{c||}{{\tt tol}$= 10^{-3}$ and $\alpha =0$} \\ \hline		
	$n$ & $q_S(M_{orig})$ & $q_S (M_{scal})$ & $\kappa (D_\ell)$ & $\kappa (D_r)$ & steps \\	\hline
 400 & 1.94e+10  & 1.18e+01 & 1.11e+04 & 1.64e+04 & 1367.5 \\
 800 & 4.90e+09 & 1.20e+01 & 5.53e+03 & 8.70e+03 & 1616 \\
 1200 & 1.12e+10 & 1.26e+01 & 7.17e+03 & 7.27e+03 & 1470.7 \\
 1600 & 2.79e+09 &  1.27e+01 & 6.14e+03 & 4.30e+03 & 1323.3  \\
 2000 & 4.07e+09 & 1.27e+01 & 5.32e+03 & 5.94e+03 & 1382.4
\end{tabular}

  } }
	\end{center}
\end{table}

\cblu{
Our second test is organized in the same way as that in Table \ref{tab:table61-fro}, but the generated matrices $A$ and $B$ are sparse, with only around 1 \% of their entries different from zero. They are generated as described in the caption of Table \ref{tab:table62-fro}. The sparsity of the corresponding $M$ matrices imply that they may have not often total support. In fact, the algorithm in Appendix A with $r = c = n \1_n$ has not converged in 2000 steps for any of the matrices $M$ generated in this test with {\tt tol}$= 10^{-3}$. This indicates that a regularization would be needed in {\em exact arithmetic} for these pencils. However, the algorithm has always converged rather quickly with {\tt tol}$= 1$, yielding, moreover, very satisfactory scalings as measured by $q_S (M_{scal})$. The results are shown in Table \ref{tab:table62-fro}. This test is just one example of a phenomenon that we have observed very often, namely, that the use of {\tt tol}$= 1$ has very often a regularization effect that makes it unnecessary to use, for computing accurate eigenvalues of pencils, the regularization techniques in Section \ref{sec:regularized}. We announced this phenomenon in Example \ref{ex-f51}, but we have observed it in many other cases where the matrix $M$ does not have total support and it has led us to make the comments in Remark \ref{rem.regular}. Observe that the convergence in Table \ref{tab:table62-fro} is slower than in Table \ref{tab:table61-fro}. As we discuss below, this is due to the fact that the values of $q_S (M_{orig})$ are larger, but also due to the larger sparsity.

}

\begin{table}[h!]
	\begin{center}
\cblu{
		\caption{Numerical test illustrating that the use of {\tt tol}$= 1$ has often a regularizing effect. The algorithm in Appendix A with $r = c = n \1_n$ has been applied to the matrices $M$ of $n\times n$ pencils generated in {\rm MATLAB} as {\tt A=eye(n) +sprandn(n,n,0.01).}\!\!\^ {\tt (20)} and {\tt B=eye(n) + sprandn(n,n,0.01).}\!\!\^{\tt (20)} with {\tt tol}$= 1$, and the results are shown in the table. In contrast, the same algorithm applied to the same set of pencils with {\tt tol}$= 10^{-3}$ does not converge in 2000 steps for any of the generated matrices. The same happens if the power $20$ is replaced by $10$.}
		\label{tab:table62-fro}
	{\small	
\begin{tabular}{c|c||c|c|c|c||} 
  & & \multicolumn{4}{c||}{{\tt tol}$= 1$ and $\alpha =0$}  \\ \hline		
	$n$ & $q_S(M_{orig})$ & $q_S (M_{scal})$ & $\kappa (D_\ell)$ & $\kappa (D_r)$ & steps  \\	\hline
400 & 2.58e+22 & 1.33e+01 & 9.19e+16 & 1.27e+17 & 40.5 \\
800 & 3.29e+24 & 1.40e+01 & 2.72e+14 & 3.74e+15 & 33.9 \\
1200 & 9.25e+25 & 1.47e+01 & 2.89e+11 & 2.97e+13 & 27.8 \\
1600 & 2.31e+26 & 1.51e+01 & 2.34e+11 & 1.11e+12 & 28.4 \\
2000 & 3.76e+22 & 1.51e+01 & 2.34e+10 & 1.61e+11 & 26
\end{tabular}
  } }
	\end{center}
\end{table}

\cblu{
Our third test is organized as the previous ones. The test pencils are in this case random permutations of square block diagonal pencils with rectangular diagonal blocks. They are generated as described in the caption of Table \ref{tab:table63-fro}. None of the corresponding $M$ matrices has in this case total support. The key difference with respect to the tests in Tables \ref{tab:table61-fro} and \ref{tab:table62-fro} is that in this case the algorithm in Appendix A with {\tt tol}$= 1$ and $r = c = n \1_n$ applied to the matrices $M$ never converges in 2000 steps, i.e., {\tt tol}$= 1$ does not have a regularizing effect for these pencils. Thus, the use of the regularization is mandatory in this case. The results are shown in Table \ref{tab:table63-fro}. We emphasize two main points on the results. First, though the obtained values for $q_S (M_{scal})$ are much better than those of $q_S (M_{orig})$, they are far from $1$. Moreover, the values of $q_S (M_{scal})$ do not improve by decreasing the value of $\alpha$. Despite these facts, we will see in some experiments done in  Subsection \ref{sub:3} on similar pencils, that the regularized scaling has significant positive effects on the accuracy of the computed eigenvalues.}

\begin{table}[h!]
	\begin{center}
\cblu{
		\caption{Numerical test illustrating pencils where the regularization is mandatory even if {\tt tol}$= 1$ is used. The algorithm in Appendix A with $r = c = n \1_n$ applied to the matrices $M$ of $n\times n$ pencils generated in {\rm MATLAB} as random permutations of {\tt A = blkdiag(randn(n1,n2).}\!\!\^ {\tt (20) , randn(n2,n1).}\!\!\^ {\tt (20))} and {\tt B = blkdiag(randn(n1,n2).}\!\!\^ {\tt (20) , randn(n2,n1).}\!\!\^ {\tt (20))} with {\tt n1 = n/5} and {\tt n2 = n - n1} does not converge in 2000 steps with {\tt tol}$= 1$. The same happens if the exponent 20 is replaced by 10 or 5. In contrast, the algorithm in Appendix A with $r = c = (2 n) \1_{2 n}$ applied to the matrices $M_\alpha^{\circ 2}$ in \eqref{eq:Malpha} with {\tt tol}$= 1$ and $\alpha = 0.5, 10^{-4}$ converges and the results are shown below. We have checked that the use of smaller values of $\alpha$ does not improve the quality of the achieved scaling, but worsens the condition numbers of $D_\ell$ and $D_r$ and increases the number of steps until convergence}
		\label{tab:table63-fro}
	{\small	
\begin{tabular}{c|c||c|c|c|c||} 
  & & \multicolumn{4}{c||}{{\tt tol}$= 1$ and $\alpha =0.5$}  \\ \hline		
	$n$ & $q_S(M_{orig})$ & $q_S (M_{scal})$ & $\kappa (D_\ell)$ & $\kappa (D_r)$ & steps  \\	\hline
200 & 1.59e+15 & 1.75e+09 & 3.74e+12 & 3.30e+12 & 32.8 \\
400 & 2.14e+14 & 4.93e+07 & 1.94e+13 & 6.77e+13 & 34.5 \\
600 & 9.30e+13 & 6.64e+06 & 1.37e+14 & 5.63e+13 & 33.6 \\
800 & 1.20e+13 & 3.45e+06 & 1.20e+14 & 1.13e+14 & 33.0 \\
1000 & 4.72e+12 & 4.57e+06 & 1.48e+14 & 1.41e+14 & 34.0
\end{tabular}

\vspace*{0.3cm}

\begin{tabular}{c|c||c|c|c|c||} 
  & & \multicolumn{4}{c||}{{\tt tol}$= 1$ and $\alpha = 10^{-4}$}  \\ \hline		
	$n$ & $q_S(M_{orig})$ & $q_S (M_{scal})$ & $\kappa (D_\ell)$ & $\kappa (D_r)$ & steps  \\	\hline
200 & 1.59e+15 & 1.83e+09 & 2.25e+16 & 1.46e+16 & 39.5 \\
400 & 2.14e+14 & 5.42e+07 & 9.01e+16 & 5.22e+17 & 41.0 \\
600 & 9.30e+13 & 4.97e+06 & 6.20e+17 & 2.59e+17 & 40.0 \\
800 & 1.20e+13 & 3.72e+06 & 5.04e+17 & 4.76e+17 & 39.1 \\
1000 & 4.72e+12 & 5.20e+06 & 7.21e+17 & 5.76e+17 & 40.1
\end{tabular}

  } }
	\end{center}
\end{table}

\cblu{
We finish this subsection with two additional tests. The first one is described and reported in Table \ref{tab:table64-fro} and is as the one in Table \ref{tab:table61-fro} but with starting matrices $M$ that are less strongly unbalanced as measured by $q_S(M_{orig})$. This leads to a much faster  convergence than in Table \ref{tab:table61-fro}, as it is naturally expected. The comparison of Tables \ref{tab:table61-fro} and \ref{tab:table64-fro} shows that the number of steps until converges grows with the unbalancing of the $M$ matrices but, also, that is independent of the dimension of the matrices. The last test is described and reported in Table \ref{tab:table65-fro} and is as the one in Table \ref{tab:table62-fro} but with sparse starting matrices $M$ that are less strongly unbalanced, which lead again to a much faster convergence, independent, more or less, of the dimension of the matrices. The comparison of Table \ref{tab:table61-fro} (for {\tt tol}$= 1$), for dense pencils, and of Table \ref{tab:table65-fro}, for sparse pencils, is interesting because both show similar values of $q_S(M_{orig})$ but the convergence is slower in the sparse case. This illustrates that for {\tt tol}$= 1$, the well-known effect that sparsity slows down the convergence of the Sinkhorn-Knopp algorithm also holds \cite{Knight}.}

\begin{table}[h!]
	\begin{center}
\cblu{
		\caption{Numerical test equal to that in Table \ref{tab:table61-fro} for {\tt tol}$= 1$ except for the fact that the $n\times n$ pencils are generated in {\rm MATLAB} as {\tt A=randn(n,n).}\!\!\^ {\tt (10)} and {\tt B=randn(n,n).}\!\!\^{\tt (10)}. The use of the exponent 10 instead of 20 in the generation of the test matrices implies that the original matrices $M$ are better equilibrated than those in Table \ref{tab:table61-fro}, as indicated by the values of $q_S(M_{orig})$, which, in turns, implies a faster convergence in approximately half of the steps}
		\label{tab:table64-fro}
	{\small	
\begin{tabular}{c|c||c|c|c|c||} 
  & & \multicolumn{4}{c||}{{\tt tol}$= 1$ and $\alpha =0$}  \\ \hline		
	$n$ & $q_S(M_{orig})$ & $q_S (M_{scal})$ & $\kappa (D_\ell)$ & $\kappa (D_r)$ & steps  \\	\hline
 400 & 4.04e+04 & 1.11e+01 & 3.84e+01 & 8.64e+01 & 5.1 \\
 800 & 1.75e+04 & 1.07e+01 & 2.72e+01 & 6.40e+01 & 5.0 \\
 1200 & 1.78e+04 & 1.12e+01 & 2.24e+01 & 5.12e+01 & 5.1 \\
 1600 & 1.37e+04 & 1.13e+01 & 2.72e+01 & 4.80e+01 & 4.8 \\
 2000 & 1.42e+04 & 1.15e+01 & 2.40e+01 & 5.12e+01 & 5.0
\end{tabular}
} }
	\end{center}
\end{table}

\begin{table}[h!]
	\begin{center}
\cblu{
		\caption{Numerical test equal to that in Table \ref{tab:table62-fro} except for the fact that the $n\times n$ pencils are generated in {\rm MATLAB} as {\tt A=eye(n) + sprandn(n,n,0.01).}\!\!\^ {\tt (10)} and {\tt B=eye(n) + sprandn(n,n,0.01).}\!\!\^{\tt (10)}. The use of the exponent 10 instead of 20 in the generation of the test matrices implies that the original matrices $M$ are better equilibrated than those in Table \ref{tab:table62-fro}, as indicated by the values of $q_S(M_{orig})$, which, in turns, implies a faster convergence in approximately half of the steps}
		\label{tab:table65-fro}
	{\small	
\begin{tabular}{c|c||c|c|c|c||} 
  & & \multicolumn{4}{c||}{{\tt tol}$= 1$ and $\alpha =0$}  \\ \hline		
	$n$ & $q_S(M_{orig})$ & $q_S (M_{scal})$ & $\kappa (D_\ell)$ & $\kappa (D_r)$ & steps  \\	\hline
400 & 8.87e+10 & 1.31e+01 & 1.96e+08 & 1.85e+08 & 20.6 \\
800 & 1.04e+12 & 1.38e+01 & 8.07e+06 & 1.93e+07 & 17.1 \\
1200 & 3.96e+12 & 1.41e+01 & 3.87e+05 & 3.04e+06 & 14.3 \\
1600 & 3.39e+12 & 1.42e+01 & 3.28e+05 & 5.24e+05 & 14.5 \\
2000 & 1.12e+11 & 1.54e+01 & 8.19e+04 & 2.29e+05 & 13.1
\end{tabular}
  } }
	\end{center}
\end{table}

\cblu{As a summary of the results in this subsection, we emphasize that, even for pencils leading to extremely unbalanced matrices $M$, the computational cost of the scaling procedures proposed in this paper with the stopping criterion {\tt tol}$= 1$ is much smaller than the cost of computing the eigenvalues. For brevity, results on rectangular pencils are delayed until Section \ref{sub:secrectangularnum}.}

\cblu{

\subsection{Examples on the accuracy of computed eigenvalues of regular pencils}  \label{sub:2}
In this section, we discuss numerical tests for three families of regular pencils. In each of these families, we generated random diagonalizable $n \times n$ regular pencils $\lambda B- A$ for which their ``exact'' eigenvalues $\lambda_i$ were known. Then, we applied the $QZ$-algorithm \cite{QZ} in MATLAB to such pencils, to the scaled pencils $D_\ell(\la B-A)D_r$ obtained by applying the algorithm in Appendix A with $r=c= n \1_n$ and {\tt tol}$=1$ to $M =  |A|^{\circ 2}+|B|^{\circ 2}$, and to the pencils balanced by Ward's method \cite{Ward}. In all cases, we constrained the diagonal elements of the diagonal scaling matrices to be integer powers of two. Since MATLAB does not have a built-in function implementing Ward's method, we used the one in \cite{Weder}. For each generated pencil, we compared the ``exact'' eigenvalues $\lambda_i$ of the pencil with the eigenvalues $\tilde \lambda_i$ computed via the three options described above. For the comparison of the eigenvalues,
we used their chordal distances \cite{Stewart-Sun}
$$c_i  := \chi(\la_i,\tilde \lambda_i) := \frac{|\lambda_i-\tilde\lambda_i|}{\sqrt{1+|\lambda_i|^2}\sqrt{1+|\tilde\lambda_i|^2}}.$$
We compared  the quantities $c :=\|[c_1,\ldots, c_n]\|_2$ for the
original pencil $(\la B-A)$ ($c_{orig}$), for the balanced pencil $D_\ell(\la B-A)D_r$ constructed by applying the algorithm in Appendix A with $r=c= n \1_n$ and {\tt tol}$=1$ to $M$  $(c_{bal}$) and for the balanced pencil constructed by Ward's method $(c_{ward}$). The regularization techniques of Section \ref{sec:regularized} were not used in this section since the algorithm in Appendix A applied to $M$ with {\tt tol}$=1$ always converged in a very small number of steps, as can be seen in the tables of this subsection. In fact, we have not found any {\em regular} pencil where the algorithm in Appendix A with $r=c= n \1_n$ and {\tt tol}$=1$ applied to $M$ does not converge in a small number of steps, even considering very sparse regular pencils.}

\cblu{In the first family of tests of this subsection, we generated $500 \times 500$ random diagonalizable pencils of the form $T_\ell(\lambda \Lambda_B- \Lambda_A)T_r$ where $(\lambda \Lambda_B- \Lambda_A)$ is in standard normal form \cite{LemVD}, i.e., $\Lambda_A$ and $\Lambda_B$ are diagonal, and $|\Lambda_A|^2 + |\Lambda_B|^2 = I_n$. The condition number of the random square nonsingular matrices $T_\ell$ and $T_r$ was controlled by taking the $k$th power of normally distributed random numbers $r_{ij}$ as their elements. A larger power $k$ then typically yields a larger condition number. The obtained results are shown in Table \ref{tab:table66-fro}, where each row corresponds to a value of $k$ taken in increasing order from $k = 1:5:41$ in MATLAB notation.
This experiment shows that the scaling proposed in Section \ref{sec:regular} based on the algorithm in Appendix A does improve the accuracy of the computed eigenvalues with respect to the original pencil and to the pencil scaled by Ward's method, especially when
the pencil corresponds to badly conditioned left and right diagonalizing transformations $T_\ell$ and $T_r$. Moreover, we see that the algorithm in Appendix A converged in a very small number of steps and produced a very well scaled matrix $\widetilde{M}$.}

\begin{table}[h!]
	\begin{center}
\cblu{
		\caption{Eigenvalue accuracy of the $QZ$-algorithm for regular $500 \times 500$ pencils: for the original pencil, for the pencil balanced by applying the algorithm in Appendix A with $r=c= n \1_n$ and {\tt tol}$=1$ to $M =  |A|^{\circ 2}+|B|^{\circ 2}$, and for the pencil balanced by Ward's method. The improvement in the scaling of $M$ produced by the algorithm in Appendix A is also shown in terms of $q_S(M_{orig})$ and $q_S (M_{scal})$ (see \eqref{eq.qs}), as well as the number of its steps until convergence}
		\label{tab:table66-fro}
	{\small	\begin{tabular}{c|c|c|c|c|c|c} 
$\kappa(T_\ell)$ & $\kappa(T_r)$ & $c_{orig}$ & $c_{bal}$ & $c_{ward}$ & $c_{bal}/c_{orig}$ & $c_{bal}/c_{ward}$  \\
			\hline
    2.45e+03 & 1.03e+03 & 7.42e-13 & 7.42e-13 & 7.19e-13 & 1.00e+00 & 1.03e+00 \\
       4.11e+03 & 4.20e+03 & 5.29e-13 & 4.25e-13 & 4.61e-13 & 8.03e-01 & 9.22e-01 \\
       2.01e+05 & 5.26e+04 & 1.59e-11 & 4.89e-12 & 5.33e-12 & 3.08e-01 & 9.17e-01 \\
       4.25e+07 & 3.87e+06 & 9.94e-10 & 1.92e-11 & 2.28e-10 & 1.93e-02 & 8.39e-02 \\
       4.55e+08 & 2.83e+07 & 2.09e-08 & 1.07e-10 & 2.07e-09 & 5.13e-03 & 5.20e-02 \\
       7.47e+10 & 2.62e+10 & 1.19e-05 & 5.67e-08 & 8.97e-06 & 4.76e-03 & 6.31e-03 \\
       9.18e+11 & 7.91e+11 & 2.57e-03 & 1.96e-05 & 1.16e-03 & 7.63e-03 & 1.69e-02 \\
       5.31e+14 & 1.29e+14 & 4.80e-01 & 3.44e-06 & 7.40e-03 & 7.18e-06 & 4.65e-04 \\
       9.66e+16 & 5.23e+14 & 1.33e-01 & 2.20e-03 & 2.09e-01 & 1.65e-02 & 1.05e-02
		\end{tabular}

\vspace*{0.3cm}

\begin{tabular}{c|c|c} 
$q_S(M_{orig})$ & $q_S (M_{scal})$ & steps \\ \hline
1.62e+00 & 1.62e+00 & 1 \\
       4.25e+03 & 5.73e+00 & 4 \\
       1.11e+06 & 9.03e+00 & 5 \\
       9.32e+09 & 1.16e+01 & 9 \\
       6.12e+11 & 1.01e+01 & 12 \\
       7.54e+16 & 9.97e+00 & 14 \\
       5.57e+18 & 1.18e+01 & 17 \\
       5.15e+24 & 1.07e+01 & 23 \\
       3.53e+26 & 1.35e+01 & 21
\end{tabular}}

}
	\end{center}
\end{table}

\cblu{It is well known that Ward's method can severely deteriorate the accuracy of the computed eigenvalues of some pencils \cite[Ch. 2, Sect. 4.2]{kressner}, \cite{LemVD}. In the second family of tests of this subsection, we generated a family of $500 \times 500$ pencils where Ward's method led to computed eigenvalues with large errors but the method in Section \ref{sec:regular} performed very well in accuracy and convergence rate. We emphasize that we have not been able to generate pencils with the opposite behavior. The pencils were generated as follows: (1) a random $500 \times 500$ matrix $T$ was constructed with the MATLAB command {\tt randn}; (2) small entries were created in $T$ with $T(1,2:500) = 10^{-k} T(1,2:500)$ and $T(4:500,3) = 10^{-k} T(4:500,3)$; (3) take $A = T D$, with $D$ a random diagonal matrix of integer positive numbers, and $B = T$. Observe that the eigenvalues of $\la B - A$ are precisely the diagonal entries of $D$. The results are shown in Table \ref{tab:table67-fro}, where each row corresponds to a value of $k$ taken from $k = 1:2:11$ in MATLAB notation.
}

\begin{table}[h!]
	\begin{center}
\cblu{
		\caption{Eigenvalue accuracy of the $QZ$-algorithm for regular $500 \times 500$ pencils for which Ward's method deteriorates the precision of computed eigenvalues: for the original pencil, for the pencil balanced by applying the algorithm in Appendix A with $r=c= n \1_n$ and {\tt tol}$=1$ to $M$, and for the pencil balanced by Ward's method}
		\label{tab:table67-fro}
	{\small	\begin{tabular}{c|c|c|c|c|c|c} 
k & $c_{orig}$ & $c_{bal}$ & $c_{ward}$ & $c_{bal}/c_{orig}$ & $c_{bal}/c_{ward}$ &
  $c_{ward}/c_{orig}$\\
			\hline
     1 & 2.61e-13 & 3.40e-15 & 8.87e-15 & 1.31e-02 & 3.84e-01 & 3.40e-02 \\
       3 & 1.48e-13 & 7.59e-15 & 1.91e-14 & 5.14e-02 & 3.98e-01 & 1.29e-01 \\
       5 & 4.13e-13 & 8.72e-15 & 4.56e-09 & 2.11e-02 & 1.91e-06 & 1.10e+04 \\
       7 & 7.16e-14 & 2.27e-15 & 3.47e-02 & 3.17e-02 & 6.54e-14 & 4.84e+11 \\
       9 & 3.90e-13 & 3.01e-15 & 1.05e+00 & 7.72e-03 & 2.87e-15 & 2.69e+12 \\
       11 & 1.34e-13 & 7.99e-15 & 1.08e+00 & 5.96e-02 & 7.38e-15 & 8.08e+12
		\end{tabular}

\vspace*{0.3cm}

\begin{tabular}{c|c|c} 
$q_S(M_{orig})$ & $q_S (M_{scal})$ & steps \\ \hline
5.11e+04 & 4.02e+00 & 2 \\
       1.16e+05 & 4.33e+00 & 3 \\
       1.43e+05 & 4.33e+00 & 3 \\
       6.40e+03 & 4.50e+00 & 3 \\
       1.47e+05 & 4.73e+00 & 3 \\
       1.37e+05 & 4.74e+00 & 3
\end{tabular}}

}
	\end{center}
\end{table}

\cblu{In the experiments presented so far in this subsection, the scaling method in Section \ref{sec:regular} always improved significantly the accuracy of the computed eigenvalues with respect to the original unscaled pencil. However, there are pencils where the improvement is much larger. This is illustrated in the last family of tests of this subsection. The pencils were constructed as those in the experiment of Table \ref{tab:table66-fro}, i.e., $T_\ell(\lambda \Lambda_B- \Lambda_A)T_r$, but with different $T_\ell$ and $T_r$. In this case, $T_\ell = D_1 Q_\ell$ and $T_r = Q_r D_2$, with $Q_\ell$ and $Q_r$ random orthogonal matrices and $D_1$ and $D_2$ random diagonal matrices with condition numbers $10^k$ and geometrically distributed singular values, constructed with the command {\tt gallery('randsvd',...)} of MATLAB. The results are shown in Table \ref{tab:table68-fro} for $1000 \times 1000$ pencils and $k = 1, 10, 19$ (each value for each row of the table). Ward's method also yields very accurate eigenvalues.
}

\begin{table}[h!]
	\begin{center}
\cblu{
		\caption{Eigenvalue accuracy of the $QZ$-algorithm for regular $1000 \times 1000$ pencils for which the method based on the algorithm in Appendix A applied to $M$ and Ward's method work both very well}
		\label{tab:table68-fro}
	{\small	\begin{tabular}{c|c|c|c|c|c|c|c} 
  $c_{orig}$ & $c_{bal}$ & $c_{ward}$ & $\!\! c_{bal}/c_{orig}\!\!$ & $\!\! c_{bal}/c_{ward} \!\!$ & $q_S(M_{orig})$ & $\!\! q_S (M_{scal}) \!\!$ & \!\! steps \!\! \\
			\hline
       5.27e-14 & 5.33e-14 & 5.17e-14 & 1.01e+00 & 1.03e+00 & 1.12e+02 & 4.13e+00 & 2\\
       4.47e-06 & 5.23e-14 & 5.77e-14 & 1.17e-08 & 9.05e-01 & 1.65e+20 & 4.21e+00 & 3  \\
        1.33e+01 & 6.49e-14 & 6.41e-14 & 4.86e-15 & 1.01e+00 & 1.53e+38 & 4.13e+00 & 3
		\end{tabular}}
}
	\end{center}
\end{table}

\cblu{As a consequence of the results in this subsection, we emphasize again that the scaling method in Section \ref{sec:regular}, i.e., that in \cite{LemVD}, often contributes to improve the accuracy of computed eigenvalues of regular pencils significantly and outperforms the method of Ward \cite{Ward}, which is the only one available so far in LAPACK \cite{lapack}.}

\cblu{
\subsection{Examples on the accuracy of computed eigenvalues of singular square pencils}  \label{sub:3}
In this section, we discuss tests for two families of singular square pencils. The first family includes dense pencils for which the regularization in Section \ref{sec:regularized} is not needed, while the second one corresponds to sparse pencils for which the regularization is necessary. For completeness, Ward's method is also considered in the comparisons, because, although it was developed for regular pencils, it has worked on the singular ones of this subsection. As in Subsection \ref{sub:2}, we generated random singular pencils whose ``exact'' eigenvalues are known and we used the vectors of chordal distances, $c :=\|[c_1,\ldots, c_n]\|_2$ for the
original pencil $(\la B-A)$ ($c_{orig}$), for the balanced pencil $D_\ell(\la B-A)D_r$ constructed by the methods in either Section \ref{sec:regular} or \ref{sec:regularized}, and for the balanced pencil constructed by Ward's method $(c_{ward}$), in order to check the improvements that the different scalings produced on the accuracy of the computed eigenvalues.}

\cblu{The first family of dense pencils considered in this subsection is constructed in the same way as the pencils in Table \ref{tab:table66-fro}, but} we replaced one of the diagonal pairs of the \cblu{$500 \times 500$} pencil $(\lambda \Lambda_B- \Lambda_A)$ generated in the regular example
by two zeros, \cblu{thus creating} a singular pencil. Each transformed pencil $(\la B -A) := T_\ell(\lambda \Lambda_B- \Lambda_A)T_r$
is therefore also singular, but its left and right rational null spaces are both of dimension 1 and their minimal bases are formed by constant vectors \cite{Van79}. For that reason, the
regular part of that singular pencil has dimension \cblu{$499\times 499$} and its eigenvalues are the remaining \cblu{499} eigenvalues of
$(\lambda \Lambda_B- \Lambda_A)$. If we follow the same procedure as in the \cblu{regular} experiment, the $QZ$-algorithm applied to $(\la B -A)$
should in principle yield arbitrary eigenvalues, since it is known that the $QZ$-algorithm is backward stable and that there exist arbitrarily small perturbations of
square singular pencils that make them regular, but with arbitrary spectrum in the complex plane \cite{Van79}. However, it has been shown that such perturbations are very particular, and that, generically, tiny perturbations of a singular square pencil makes it regular with eigenvalues that are tiny perturbations of the eigenvalues of the unperturbed singular pencil, together with some other ``arbitrary'' eigenvalues determined by the perturbation \cite{detedop2010,detedop2008}. Even more, starting from these ideas, it has been shown very recently that it is possible to define sensible and useful ``weak'' condition numbers for the eigenvalues of a singular square pencil \cite{LotN20}. This explains the well-known fact that, in practice, the $QZ$-algorithm applied to a singular square matrix pencil finds almost always its eigenvalues, albeit with some loss of accuracy. Therefore, it makes sense to apply the $QZ$ algorithm to our generated singular pencils as well as to their scaled versions.
The numerical results are reported in \cblu{Table \ref{tab:table69-fro}, where each row corresponds to a value of $k$ taken in increasing order from $k = 1 : 5 : 41$ as in Table \ref{tab:table66-fro}}. We generated the data
just as in the \cblu{experiment for regular pencils in Table \ref{tab:table66-fro}}, except for the one eigenvalue replaced by $0/0$ or, in other words, by NaN. When comparing the ``original'' spectrum with the computed one, we excluded NaN in the original set
and looked for the best matching \cblu{499} eigenvalues in the ``computed'' spectrum. \cblu{It is clear from Table \ref{tab:table69-fro} that the balancing proposed in Section \ref{sec:regular} also improves the accuracy of the computed
eigenvalues of singular square pencils, both with respect to the original pencil and with respect to the one balanced by Ward's method, and that needs a small number of steps to converge.}

\begin{table}[h!]
	\begin{center}
\cblu{
		\caption{Eigenvalue accuracy of the plain $QZ$-algorithm for singular $500 \times 500$ dense pencils: for the original pencil, for the pencil balanced by applying the algorithm in Appendix A with $r=c= n \1_n$ and {\tt tol}$=1$ to $M =  |A|^{\circ 2}+|B|^{\circ 2}$, and for the pencil balanced by Ward's method. The improvement in the scaling of $M$ produced by the algorithm in Appendix A is also shown in terms of $q_S(M_{orig})$ and $q_S (M_{scal})$ (see \eqref{eq.qs}), as well as the number of its steps until convergence}
		\label{tab:table69-fro}
	{\small	\begin{tabular}{c|c|c|c|c|c|c} 
$\kappa(T_\ell)$ & $\kappa(T_r)$ & $c_{orig}$ & $c_{bal}$ & $c_{ward}$ & $c_{bal}/c_{orig}$ & $c_{bal}/c_{ward}$  \\
			\hline
   4.30e+03 & 4.10e+03 & 1.88e-12 & 1.88e-12 & 8.27e-12 & 1.00e+00 & 2.28e-01 \\
       1.69e+04 & 2.12e+04 & 1.77e-11 & 1.85e-12 & 6.17e-12 & 1.04e-01 & 2.99e-01 \\
       1.06e+06 & 9.83e+04 & 1.88e-11 & 1.19e-11 & 5.04e-12 & 6.34e-01 & 2.37e+00 \\
       7.47e+05 & 2.73e+06 & 1.98e-10 & 1.40e-10 & 7.13e-11 & 7.08e-01 & 1.97e+00 \\
       1.20e+08 & 6.49e+08 & 1.62e-08 & 4.13e-11 & 4.13e-09 & 2.55e-03 & 9.99e-03 \\
       2.32e+10 & 2.75e+09 & 5.20e-07 & 5.00e-09 & 2.15e-07 & 9.62e-03 & 2.33e-02 \\
       3.59e+13 & 2.59e+12 & 3.25e-03 & 2.83e-07 & 5.40e-05 & 8.71e-05 & 5.24e-03 \\
       1.63e+16 & 3.03e+13 & 3.46e-02 & 3.55e-05 & 3.84e-03 & 1.03e-03 & 9.25e-03 \\
       1.63e+18 & 1.48e+14 & 8.15e-02 & 9.12e-06 & 1.22e-02 & 1.12e-04 & 7.46e-04
		\end{tabular}

\vspace*{0.3cm}

\begin{tabular}{c|c|c} 
$q_S(M_{orig})$ & $q_S (M_{scal})$ & steps \\ \hline
 1.57e+00 & 1.57e+00 & 1 \\
       1.09e+03 & 6.04e+00 & 3 \\
       1.40e+06 & 8.91e+00 & 7 \\
       2.66e+09 & 9.43e+00 & 8 \\
       1.10e+12 & 1.01e+01 & 13 \\
       1.31e+16 & 9.26e+00 & 13 \\
       1.13e+20 & 1.43e+01 & 16 \\
       1.72e+25 & 1.20e+01 & 17 \\
       2.11e+26 & 1.16e+01 & 18
\end{tabular}}

}
	\end{center}
\end{table}

Though the direct use of the $QZ$-algorithm is a simple option for computing the eigenvalues of a singular square pencil when the accuracy requirements are moderate, the correct handling of a singular pencil is to first ``deflate'' its left and right null spaces, and
then compute the spectrum of the regular part of that singular pencil, i.e., to apply the staircase algorithm (see \cite{Van79}). In this experiment,
it turns out that the left and right null spaces are one-dimensional and are given, respectively, by the left null
vector of $\left[\begin{smallmatrix} A & B \end{smallmatrix}\right]$, and by the right null vector of
$\left[\begin{smallmatrix} A \\ B \end{smallmatrix}\right]$, which we both computed using a singular value
decomposition of these compound matrices. \cblu{After this deflation was applied to the original pencil $(\la B-A)$, to the
pencil $D_\ell(\la B-A)D_r$ scaled by the method in Section \ref{sec:regular} and to the one balanced by Ward's method,} we again computed the spectrum of the deflated pencils with the $QZ$-algorithm. The results for the same
data as reported in Table \ref{tab:table69-fro} are now reported in Table \ref{tab:table610-fro}. \cblu{The results in this case are similar in both tables.
We also added three columns with the sensitivities of the deflation in the original pencil $\gamma_{orig}$ and in the balanced ones by the method in Section \ref{sec:regular} and Ward's method, $\gamma_{bal}$ and $\gamma_{ward}$.} We measured the sensitivity of the left and right null vectors defining the deflation of a singular pencil $\lambda B-A$,
by
\begin{equation} \label{eq.gammas}
 \gamma := \max(\frac{\sigma_{n}\left[\begin{smallmatrix} A \\ B \end{smallmatrix}\right]}{
	\sigma_{n-1}\left[\begin{smallmatrix} A \\ B \end{smallmatrix}\right]},
\frac{\sigma_{n}\left[\begin{smallmatrix} A & B \end{smallmatrix}\right]}{
	\sigma_{n-1}\left[\begin{smallmatrix} A & B \end{smallmatrix}\right]}),
\end{equation}
i.e. the largest ratio between the two smallest singular values of the matrices that define these null vectors.
It is an indication about how much these vectors can rotate when perturbing the pencil. It is easy to see from
\cblu{the data that the accuracy of the computed eigenvalues of the deflated pencil is closely related to the sensitivity of the deflation itself.}

\begin{table}[h!]
	\begin{center}
\cblu{
		\caption{Eigenvalue accuracy of the staircase algorithm for exactly the same singular $500 \times 500$ dense pencils of Table \ref{tab:table69-fro}}
		\label{tab:table610-fro}
	{\small	\begin{tabular}{c|c|c|c|l|c|c|c} 
 $c_{orig}$ & $c_{bal}$ & $c_{ward}$ & $c_{bal}/c_{orig}$ & $\!c_{bal}/c_{ward}\!$ & $\gamma_{orig}$ &  $\gamma_{bal}$ & $\gamma_{ward}$   \\
			\hline
   2.23e-13 & 2.23e-13 & 2.33e-13 & 1.0e+00 & 9.57e-01 & 5.79e-13 & 5.79e-13 & 6.59e-13 \!\!\\
       4.53e-13 & 4.68e-13 & 2.89e-13 & 1.03e+00 & 1.62e+00 & 1.48e-11 & 4.96e-12 & 7.49e-12 \!\! \\
       6.92e-13 & 9.11e-13 & 2.00e-12 & 1.32e+00 & 4.56e-01 & 2.37e-10 & 2.33e-12 & 1.17e-10 \!\! \\
       8.36e-11 & 1.63e-11 & 9.76e-12 & 1.95e-01 & 1.67e+00 & 1.33e-07 & 5.10e-11 & 4.84e-08 \!\! \\
       6.49e-10 & 1.12e-11 & 8.46e-11 & 1.73e-02 & 1.33e-01 & 1.24e-06 & 2.02e-11 & 1.17e-07 \!\! \\
       1.41e-07 & 5.06e-09 & 2.03e-07 & 3.59e-02 & 2.49e-02 & 1.53e-03 & 2.02e-09 & 1.35e-04 \!\! \\
       7.22e-04 & 1.42e-06 & 1.03e-06 & 1.96e-03 & 1.38e+00 & 9.62e-01 & 3.43e-07 & 1.99e-01 \!\! \\
       3.33e-02 & 1.25e-06 & 9.17e-03 & 3.76e-05 & 1.36e-04 & 2.51e-01 & 2.18e-06 & 5.24e-01 \\
       1.08e-01 & 4.31e-07 & 1.84e-03 & 3.97e-06 & 2.34e-04 & 3.87e-01 & 4.59e-07 & 7.27e-01 \!\!
		\end{tabular}
}
}
	\end{center}
\end{table}

\cblu{The second family of sparse singular pencils considered in this subsection is a family of $400 \times 400$ permuted block diagonal pencils generated as follows. Set, for simplicity, $m_1 = 140$ and $n_1 = 260$. Then
\begin{equation} \label{eq.expersingsquaresparse}
\la B - A := P \begin{bmatrix}
                \la B_1 - A_1 &  \\
                 & \la B_2 - A_2
              \end{bmatrix} Q,
\end{equation}
with $P, Q$ random $400 \times 400$ permutation matrices and
\begin{eqnarray*}
\la B_1 - A_1 & = & T_{\ell 1} \begin{bmatrix}
                \la \Lambda_{B 1} - \Lambda_{A 1} &  \\
                 & 0_{1 \times (n_1 - m_1 +1)}
              \end{bmatrix} T_{r 1}, \\
\la B_2 - A_2 & = & T_{\ell 2} \begin{bmatrix}
                \la \Lambda_{B 2} - \Lambda_{A 2} &  \\
                 & 0_{(n_1 - m_1 +1) \times 1}
              \end{bmatrix} T_{r 2},
\end{eqnarray*}
where $\la \Lambda_{B 1} - \Lambda_{A 1}, \la \Lambda_{B 2} - \Lambda_{A 2}$ are random $(m_1 - 1) \times (m_1 -1)$ diagonal regular pencils in standard normal form \cite{LemVD} which contain the ``exact'' eigenvalues of $\la B - A$, and the entries of $T_{\ell 1} \in \mathbb{R}^{m_1 \times m_1}, T_{r 2} \in \mathbb{R}^{m_1 \times m_1}, T_{\ell 2} \in \mathbb{R}^{n_1 \times n_1}, T_{r 1} \in \mathbb{R}^{n_1 \times n_1}$ are
$k$th powers of normally distributed random numbers, for $k = 1 : 5 : 41$. Observe that the normal rank \cite{Van79} of these pencils is $rg = 2(m_1 -1) = 278$, that their left and right rational null spaces are both of dimension $122$ and that their minimal bases are formed by constant vectors. This mean that they are given again, respectively, by the left null
vectors of $\left[\begin{smallmatrix} A & B \end{smallmatrix}\right]$, and by the right null vectors of
$\left[\begin{smallmatrix} A \\ B \end{smallmatrix}\right]$, which were computed again using a singular value decomposition of these compound matrices. This allowed us to deflate these right and left null spaces and to obtain the regular parts of such pencils by multiplying $\la B-A$ on the left by the $rg$ left singular vectors of $\left[\begin{smallmatrix} A & B \end{smallmatrix}\right]$ corresponding to its $rg$ largest singular values and on the right by the $rg$ right singular vectors of $\left[\begin{smallmatrix} A \\ B \end{smallmatrix}\right]$ corresponding to its $rg$ largest singular values. The application of the $QZ$ algorithm to these regular parts yielded the eigenvalues of these highly singular pencils and we did it for the original pencil $(\la B-A)$, for the pencil $D_\ell(\la B-A)D_r$ scaled by the {\em regularized} method in Section \ref{sec:regularized} and for the one balanced by Ward's method.  The plain $QZ$ algorithm can also be applied directly to the pencils in \eqref{eq.expersingsquaresparse}, but it produces much larger errors than the staircase algorithm described above due to the high singularity of these pencils. The results for the staircase algorithm are shown in Table \ref{tab:table611-fro}, where each row corresponds to a value of $k$, and are discussed in the next paragraph.}

\begin{table}[h!]
	\begin{center}
\cblu{
		\caption{Eigenvalue accuracy of the staircase algorithm for singular $400 \times 400$ sparse pencils:  for the original pencil, for the pencil balanced by applying the algorithm in Appendix A with $r=c= (2n) \1_{2n}$ and {\tt tol}$=1$ to $M_\alpha^{\circ 2}$ in \eqref{eq:Malpha} with $\alpha = 0.5$, and for the pencil balanced by Ward's method. The improvement in the scaling of $M =  |A|^{\circ 2}+|B|^{\circ 2}$ produced by the algorithm in Appendix A applied to $M_\alpha^{\circ 2}$ is also shown in terms of $q_S(M_{orig})$ and $q_S (M_{scal})$, as well as the number of its steps until convergence. The last column of the second table shows that the plain $QZ$-algorithm produces much larger errors for these pencils. For brevity this is only shown for the pencils balanced by the algorithm in Appendix A, but the same happens for the other pencils}
		\label{tab:table611-fro}
	{\small	\begin{tabular}{c|c|c|c|l|c|c|c} 
 $c_{orig}$ & $c_{bal}$ & $c_{ward}$ & $c_{bal}/c_{orig}$ & $\!c_{bal}/c_{ward}\!$ & $\!\gamma_{orig}\!$ &  $\!\gamma_{bal}\!$ & $\!\gamma_{ward}\!$   \\
			\hline
   1.98e-14 & 2.25e-14 & 2.15e-14 & 1.14e+00 & 1.05e+00 & 1.10e-13 & 1.27e-13 & 9.51e-14 \\
       3.13e-14 & 2.10e-14 & 2.39e-14 & 6.71e-01 & 8.80e-01 & 4.29e-12 & 3.29e-13 & 1.38e-12 \\
       3.40e-12 & 4.49e-14 & 2.72e-13 & 1.32e-02 & 1.65e-01 & 3.80e-10 & 1.28e-12 & 5.64e-11 \\
       1.76e-11 & 4.76e-13 & 2.69e-12 & 2.71e-02 & 1.77e-01 & 3.70e-07 & 5.02e-11 & 1.01e-07 \\
       3.17e-08 & 9.47e-13 & 4.79e-11 & 2.99e-05 & 1.98e-02 & 2.26e-04 & 2.52e-10 & 1.90e-07 \\
       7.84e-03 & 7.43e-11 & 1.10e-08 & 9.48e-09 & 6.74e-03 & 1.0e+00 & 1.20e-09 & 1.11e-03 \\
       2.31e-04 & 1.21e-10 & 5.74e-07 & 5.23e-07 & 2.11e-04 & 1.0e+00 & 5.42e-08 & 1.34e-02 \\
       1.93e-02 & 3.32e-08 & 2.73e-02 & 1.72e-06 & 1.22e-06 & 1.0e+00 & 2.55e-06 & 1.0e+00 \\
       6.46e-01 & 4.64e-10 & 4.19e-03 & 7.17e-10 & 1.11e-07 & 1.0e+00 & 1.21e-07 & 1.0e+00
		\end{tabular}
\vspace*{0.3cm}

\begin{tabular}{c|c|c|c} 
$q_S(M_{orig})$ & $q_S (M_{scal})$ & steps & $c_{bal}$ plain $QZ$ \\ \hline
       3.99e+00 & 9.08e+00 & 16 & 8.56e-07\\
       9.51e+04 & 8.97e+01 & 30 &  8.93e-07\\
       1.75e+09 & 2.85e+03 & 45 &  8.40e-07\\
       5.84e+13 & 1.82e+05 & 66 &  7.01e-07\\
       1.69e+17 & 5.18e+05 & 81 &  2.44e-07 \\
       5.26e+23 & 1.41e+07 & 100 &  2.57e-07 \\
       1.12e+22 & 7.74e+06 & 112 & 9.03e-08 \\
       1.49e+26 & 1.74e+09 & 130 & 5.02e-02 \\
       2.75e+36 & 1.58e+12 & 149 & 5.96e-03
\end{tabular}
}
}
	\end{center}
\end{table}

\cblu{The matrices $M$ corresponding to the pencils in \eqref{eq.expersingsquaresparse} are very far from having total support and the Sinkhorn-Knopp algorithm applied to them with {\tt tol}$= 1$ did not converge because it produced diagonal matrices $D_\ell, D_r$ with zero diagonal entries due to underflows. Then, we regularized the problem by applying the algorithm in Appendix A with $r=c= (2n) \1_{2n}$ and {\tt tol}$=1$ to $M_\alpha^{\circ 2}$ in \eqref{eq:Malpha} with $\alpha = 0.5$. Observe, that this yielded factors $q_S (M_{scal})$ very far from $1$ but much smaller than the factors of the original matrices $q_S (M_{orig})$. Interestingly, the factors $q_S (M_{scal})$ did not improve by taking much smaller values of $\alpha$. Despite this fact, the impact of the regularized scaling on the accuracy of the computed eigenvalues is impressive both in comparison with the original pencils and with the pencils scaled by Ward's method. The new regularized method leads to the computation of very accurate eigenvalues in a problem which is extremely difficult in terms of the high singularity and of the high unbalancing of the considered pencils. We do not know any other method in the literature that can achieve such results. Moreover, the numbers of steps until convergence are still moderate taking into account the sparsity and the strong unbalancing of the pencils, and make the cost of the scaling considerably smaller than the cost of computing the eigenvalues. Finally note that Table \ref{tab:table611-fro} also includes the sensitivities of the deflations  $\gamma_{orig}$, $\gamma_{bal}$ and $\gamma_{ward}$ as in Table \ref{tab:table610-fro}. They were computed replacing $n-1$ and $n$ in \eqref{eq.gammas} by $rg$ and $rg + 1$, respectively, where $rg = 278$ is the normal rank of the pencils. We also observe in Table \ref{tab:table611-fro} a strong relation between the errors in the eigenvalues and the deflation sensitivities.
}

\cblu{The experiments in this section show that the balancing procedures of this paper improve the accuracy of the eigenvalue computation of square singular pencils as well as the sensitivity of the deflation of the regular part of a singular pencil}. We briefly mention that recently an alternative robust method to the staircase algorithm has been proposed for computing the eigenvalues of singular pencils \cite{hochstenbach2019}. This new method is related to the ideas in \cite{detedop2010,detedop2008,LotN20} and its accuracy will also improve by using our scaling strategies.

\cblu{
\subsection{Examples on the accuracy of computed eigenvalues of rectangular pencils} \label{sub:secrectangularnum} In this section we discuss briefly tests for two families of rectangular pencils that are related to the families in Subsection \ref{sub:3}. The first family includes dense pencils for which the regularization in Section \ref{sec:regularized} is not needed, while the second one corresponds to sparse pencils for which the regularization is necessary. Ward's method is not considered since it does not work for rectangular pencils. All the considered pencils $\la B -A$ have the minimal bases of their left and right null spaces formed by constant vectors. Thus, the computation of their eigenvalues is performed via the variant of the staircase algorithm described in the previous subsection, i.e., computing first the regular parts of these pencils with the singular value decompositions of the compound matrices $\left[\begin{smallmatrix} A & B \end{smallmatrix}\right]$ and
$\left[\begin{smallmatrix} A \\ B \end{smallmatrix}\right]$, and then applying the $QZ$-algorithm to the regular parts. We use the same notation and test magnitudes as in Subsection \ref{sub:3}.}

\cblu{In the first family of tests of this subsection, we generated $150 \times 450$ random pencils of the form $\la B- A = T_\ell\diag(\lambda \Lambda_B- \Lambda_A, 0_{1 \times 301})T_r$, where $(\lambda \Lambda_B- \Lambda_A)$ is in standard normal form, has dimension $149 \times 149$ and contains the ``exact'' eigenvalues of $\la B- A$. The elements of the random square nonsingular matrices $T_\ell \in \mathbb{R}^{150 \times 150}$ and $T_r \in \mathbb{R}^{450 \times 450}$ are $k$th powers of normally distributed random numbers for $k=1:5:41$. These pencils are dense and then the regularization in Subsection \ref{sec:regularized_rowandcolum} was not needed. The results are shown in Table \ref{tab:table612-fro} (each row corresponds to a value of $k$) and illustrate the very positive effect of the scaling technique of Section \ref{sec:nonsquare} on the accuracy of computed eigenvalues and its low computational cost.}

\begin{table}[h!]
	\begin{center}
\cblu{
		\caption{Eigenvalue accuracy of the staircase algorithm for rectangular $150 \times 450$ dense pencils:  for the original pencil and for the pencil balanced by applying the algorithm in Appendix A with $r=n \1_{m}$, $c=m \1_{n}$ and {\tt tol}$=1$ to $M =  |A|^{\circ 2}+|B|^{\circ 2}$. The improvement in the scaling produced by the algorithm in Appendix A is also shown in terms of $q_S(M_{orig})$ and $q_S (M_{scal})$, as well as the number of its steps until convergence.}
		\label{tab:table612-fro}
	{\small	\begin{tabular}{c|c|c|c|c|c|c|c} 
 $c_{orig}$ & $c_{bal}$ & $c_{bal}/c_{orig}$ & $\!\gamma_{orig}\!$ &  $\!\gamma_{bal}\!$ & $q_S(M_{orig})$ & $q_S (M_{scal})$ & steps    \\
			\hline
9.96e-15 & 9.96e-15 & 1.00e+00 & 1.01e-13 & 1.01e-13 & 2.29e+00 & 2.29e+00 & 2 \\
       1.95e-14 & 1.08e-14 & 5.52e-01 & 7.97e-13 & 1.97e-13 & 4.94e+03 & 7.77e+00 & 4 \\
       2.62e-13 & 1.06e-14 & 4.03e-02 & 3.03e-10 & 1.57e-13 & 1.22e+08 & 9.66e+00 & 7 \\
       2.27e-12 & 1.29e-14 & 5.68e-03 & 1.31e-08 & 7.73e-13 & 4.32e+11 & 1.06e+01 & 9 \\
       5.61e-09 & 1.97e-13 & 3.52e-05 & 1.39e-04 & 1.72e-11 & 1.36e+16 & 1.17e+01 & 12 \\
       1.51e-05 & 1.20e-13 & 7.97e-09 & 1.95e-01 & 5.78e-12 & 8.19e+23 & 1.07e+01 & 14 \\
       6.03e-05 & 1.08e-12 & 1.79e-08 & 8.08e-03 & 9.12e-12 & 3.51e+22 & 1.27e+01 & 21 \\
       5.49e-02 & 1.72e-11 & 3.13e-10 & 1.00e+00 & 1.36e-09 & 2.39e+29 & 1.17e+01 & 16 \\
       9.76e-02 & 8.40e-12 & 8.60e-11 & 1.00e+00 & 8.24e-10 & 1.24e+31 & 1.32e+01 & 24
		\end{tabular}
}
}
	\end{center}
\end{table}

\begin{table}[h!]
	\begin{center}
\cblu{
		\caption{Eigenvalue accuracy of the staircase algorithm for singular $700 \times 450$ sparse pencils:  for the original pencil and for the pencil balanced by applying the algorithm in Appendix A with $r=c= v$ in \eqref{eq.nonhomogv} and {\tt tol}$=1$ to $M_\alpha^{\circ 2}$ in \eqref{eq:Malpha} with $\alpha = 0.5$. The improvement in the scaling of $M =  |A|^{\circ 2}+|B|^{\circ 2}$ produced by the algorithm in Appendix A applied to $M_\alpha^{\circ 2}$ is also shown in terms of $q_S(M_{orig})$ and $q_S (M_{scal})$, as well as the number of its steps until convergence.}
		\label{tab:table613-fro}
	{\small	\begin{tabular}{c|c|c|c|c|c|c|c} 
 $c_{orig}$ & $c_{bal}$ & $c_{bal}/c_{orig}$ & $\!\gamma_{orig}\!$ &  $\!\gamma_{bal}\!$ & $q_S(M_{orig})$ & $q_S (M_{scal})$ & steps    \\
			\hline
1.43e-14 & 1.26e-14 & 8.86e-01 & 8.91e-14 & 8.70e-14 & 5.54e+01 & 3.27e+01 & 7 \\
       1.73e-14 & 1.39e-14 & 8.06e-01 & 4.05e-12 & 1.25e-13 & 4.64e+06 & 9.30e+03 & 13 \\
       2.81e-13 & 3.75e-14 & 1.34e-01 & 2.74e-10 & 1.30e-12 & 3.10e+11 & 1.80e+06 & 26 \\
       1.77e-11 & 1.98e-14 & 1.12e-03 & 3.28e-08 & 4.72e-12 & 5.14e+19 & 1.42e+10 & 32 \\
       2.42e-06 & 6.23e-14 & 2.58e-08 & 1.81e-03 & 1.27e-11 & 5.87e+28 & 1.09e+13 & 46 \\
       2.42e-02 & 1.15e-10 & 4.77e-09 & 1.00e+00 & 1.85e-08 & 4.53e+29 & 1.11e+18 & 46 \\
       1.69e-04 & 2.24e-11 & 1.32e-07 & 9.84e-01 & 1.07e-07 & 6.95e+37 & 1.42e+20 & 68 \\
       4.10e-03 & 2.83e-11 & 6.88e-09 & 1.00e+00 & 4.18e-06 & 9.32e+39 & 4.30e+22 & 84 \\
       9.91e-01 & 6.07e-11 & 6.13e-11 & 1.00e+00 & 1.03e-07 & 2.72e+44 & 9.90e+22 & 87
		\end{tabular}
}
}
	\end{center}
\end{table}

\cblu{For describing the second considered family of sparse rectangular pencils, we need the parameters $m_1 = 100, n_1 = 400, m_2 = 600$ and $n_2 = 50$. Then, the pencils have the structure of those in \eqref{eq.expersingsquaresparse} but with the following changes in $\la B_2 - A_2$: the dimension of $\la \Lambda_{B2} - \Lambda_{A2}$ becomes $(n_2 -1) \times (n_2 -1)$ and $0_{(n_1 -m_1+1) \times 1}$ is replaced by $0_{(m_2 -n_2+1) \times 1}$. This implies that $T_{\ell 2} \in \mathbb{R}^{m_2 \times m_2}$ and $T_{r 2} \in \mathbb{R}^{n_2 \times n_2}$. For these pencils the algorithm in Appendix A with $r=n \1_{m}$, $c=m \1_{n}$ and {\tt tol}$=1$ applied to $M$ did not converge and we used the scaling described in Subsection \ref{sec:regularized_rowandcolum} wit $\alpha = 0.5$. The results are shown in Table \ref{tab:table613-fro} (each row corresponds to a value of $k = 1:5:41$) and illustrate again the impressive positive effect of the new scaling technique on the accuracy of computed eigenvalues and its low computational cost. The values of $q_S (M_{scal})$ did not improve by considering very small values of $\alpha$.}

\section{Concluding remarks}  \label{sec:conclusion}
\cblu{In this paper, we developed new scaling techniques that apply to both regular and singular pencils. The techniques are based on applying the Sinkhorn-Knopp-like algorithm to certain nonnegative matrices easily constructed from the matrix coefficients of the pencil, and that depend on whether the scaling problem needs to be regularized or not. The regularization guarantees to get always a unique and bounded solution. Extensive numerical experiments confirm that the proposed techniques very often improve significantly the accuracy of computed eigenvalues of arbitrary pencils and outperform earlier methods for scaling regular pencils. Finally, the algorithms computing these scalings  have a computational cost that is much smaller than the cost of the subsequent generalized eigenvalue problem as a consequence of using in the Sinkhorn-Knopp-like algorithm a proper stopping criterion compatible with computing diagonal scalings whose diagonal entries are integer powers of $2$.}

\section*{Appendix A : Sinkhorn-Knopp-like algorithm MATLAB code with prescribed row sums and column sums} \label{appendix-A}

\begin{verbatim}
function [Md,dleft,dright,error] = rowcolsums(M,r,c,maxiter,tol)
%
% [Md,dleft,dright,error] = rowcolsums(M,r,c,maxiter,tol)
%
% implements a Sinkhorn-Knopp-like algorithm for
% scaling a non-negative mxn matrix M such that
%
%        Md:=diag(dleft)*M*diag(dright)
%
%  has column sums equal to a row vector c and
%  row sums equal to a column vector r where sum(c)=sum(r)
%
%  The iterative process is stopped as soon as the incremental
%  scalings are tol-close to the identity. The error vector
%  also shows the convergence pattern of the iterative scalings
%
%  Input : M, a nonnegative mxn matrix
%          r, a positive mx1 column vector and
%          c, a positive 1xn row vector satisfying sum(c)=sum(r)
%          maxiter, the maximum number of iterations
%          tol, a tolerance for the transformation updates
%  Output: Md, a nonnegative matrix with row sums r and column sums c
%              up to the tolerance tol
%          dleft and dright, the diagonals of the left/right scalings
%          error, the convergence error
%
[m,n]=size(M);error=[];
% First scale the matrix to have total sum(sum(M))=sum(c)=sum(r);
sumcr=sum(c);sumM=sum(sum(M));Md=M*sumcr/sumM;
dleft=ones(m,1)*sqrt(sumcr/sumM);dright=ones(1,n)*sqrt(sumcr/sumM);
% Then scale left and right to make row and column sums equal to r
% and c
for i=1:maxiter;
dr=sum(Md,1)./c;Md=Md./dr;er=min(dr)/max(dr);dright=dright./dr;
dl=sum(Md,2)./r;Md=dl.\Md;el=min(dl)/max(dl);dleft=dleft./dl;
error=[error er el];if max([1-er , 1-el]) < tol/2, break; end
end
% Finally scale the two scalings to have equal maxima
scaled=sqrt(max(dright)/max(dleft));
dleft=dleft*scaled;dright=dright'/scaled;
end
\end{verbatim}

\section*{Appendix B : Proof of Lemma \ref{fullyindecomp}}

\begin{proof}
	$M_{\alpha}^{\circ 2}$ has total support for all $\alpha\neq  0$ since every nonzero element is an element of a positive diagonal. To see that $M_{\alpha}^{\circ 2}$ is fully indecomposable, we apply \cite[Theorem 1.3.7]{Brualdi}. This theorem states that a square matrix with total support is fully indecomposable if and only if its bipartite graph is connected.  Then we consider the bipartite graph of $M_{\alpha}^{\circ 2},$ denoted by $BG(M_{\alpha}^{\circ 2}).$ We assume without lost of generality that $m_{1n}$ is a nonzero element of $M:=[m_{ij}].$ Then we consider the matrix $$N:=\left[\begin{array}{c|c} \frac{\alpha^2}{m^2}  1_m 1_m^{T}  &  \begin{array}{cc}
	0 & m_{1n} \\
	0 & 0
	\end{array}  \\ \hline
	\begin{array}{cc}
	0 & 0\\
	m_{1n}  & 0
	\end{array}  &  \frac{\alpha^2}{n^2}   1_n 1_n^{T}
	\end{array}\right].$$ Notice that $BG(N)$ is a sub-graph of $BG(M_{\alpha}^{\circ 2}).$ Moreover, if $\{u_1,u_2,\ldots,u_{m+n}\} $ and $\{v_1,v_2,\ldots,v_{m+n}\} $  are the sets of vertices associated with the rows and columns of $N,$ respectively, then $ BG(N)$ is of the form
	\begin{center}
		\includegraphics[scale=0.5]{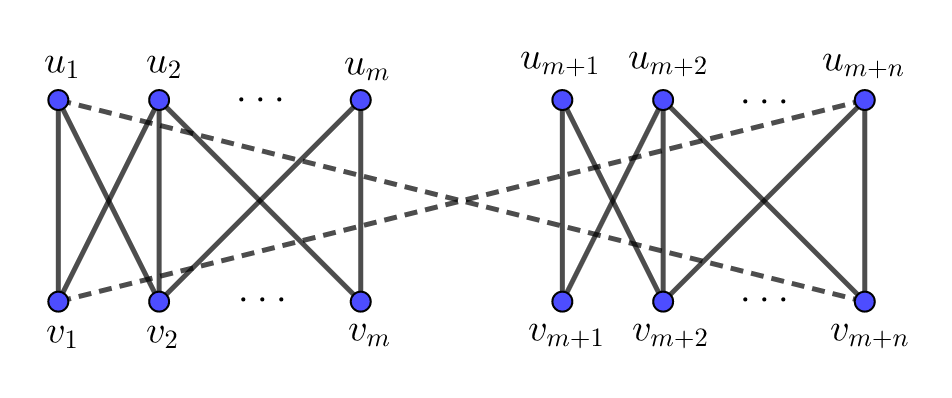}
	\end{center}	
	where the left and right groups of solid edges are each bicliques (and hence connected) and where the two dashed edges correspond to the element $m_{1n}.$ This proves that $ BG(N)$ is connected, since the dashed edges make a connection between two connected components. Therefore, $BG(M_{\alpha}^{\circ 2})$ is connected and, by \cite[Theorem 1.3.7]{Brualdi}, $M_{\alpha}^{\circ 2}$ is fully indecomposable.
\end{proof}

\vspace*{0.25cm}

\noindent \cblu{{\bf Acknowledgements.} The authors sincerely thank two anonymous referees for pointing out several significant suggestions and a number of relevant references that have contributed to improve this manuscript.}

\vspace*{0.25cm}


\end{document}